\documentclass{article}

\usepackage{moreverb,url}

\usepackage[colorlinks,bookmarksopen,bookmarksnumbered,citecolor=red,urlcolor=red]{hyperref}
\usepackage{doi}

\newcommand\BibTeX{{\rmfamily B\kern-.05em \textsc{i\kern-.025em b}\kern-.08em
T\kern-.1667em\lower.7ex\hbox{E}\kern-.125emX}}

\usepackage[algo2e,linesnumbered,vlined,algoruled]{algorithm2e}
\usepackage{geometry,xcolor,tabularx,nccmath}
\usepackage{amsmath,amsthm,amssymb,enumitem}
\usepackage{booktabs,caption,subcaption,mathtools,multirow}

\DeclareMathAlphabet{\mathdutchcal}{U}{dutchcal}{m}{n}
\SetMathAlphabet{\mathdutchcal}{bold}{U}{dutchcal}{b}{n}
\DeclareMathAlphabet{\mathdutchbcal}{U}{dutchcal}{b}{n}

\newcommand{\sE}{\mathdutchcal{E}}

\newcommand{\sG}{\mathdutchcal{G}}
\newcommand{\sV}{\mathdutchcal{V}}

\usepackage[nameinlink,capitalise]{cleveref}
\hypersetup{colorlinks={true},allcolors={blue}}

\crefformat{section}{\S#2#1#3}
\crefformat{subsection}{\S#2#1#3}
\crefformat{equation}{#2Equation~\textup{(#1)}#3}
\crefformat{cons}{#2\textbf{Constraint}~\textup{\textbf{(#1)}}#3}
\labelcrefformat{cons}{\textup{\textbf{#2(#1)#3}}}
\crefformat{consset}{#2\textbf{Constraints}~\textup{\textbf{(#1)}}#3} %
\labelcrefformat{consset}{\textup{\textbf{#2(#1)#3}}}
\crefformat{func}{#2Function~\textup{(#1)}#3}
\labelcrefformat{func}{#2Function~\textup{(#1)}#3}
\crefformat{algo}{#2\textbf{Algorithm}~\textup{\textbf{#1}}#3}
\crefformat{objfunc}{#2\textbf{Objective function}~\textup{\textbf{(#1)}}#3}
\labelcrefformat{objfunc}{\textup{(#2\textbf{#1}#3)}}
\crefformat{ch}{#2Chapter~\textup{#1}#3}
\crefformat{sec}{#2Section~\textup{#1}#3}
\crefformat{tab}{#2\textbf{Table}~\textup{\textbf{#1}}#3}
\crefformat{fig}{#2\textbf{Figure}~\textup{\textbf{#1}}#3}
\Crefformat{figure}{#2\textbf{Figure}~\textup{\textbf{#1}}#3}
\crefdefaultlabelformat{#2\textbf{#1}#3} %
\Crefname{figure}{\textbf{Figure}}{\textbf{Figures}}
\crefformat{eq}{#2equation~\textup{(#1)}#3}
\labelcrefformat{eq}{\textup{(#2#1#3)}}
\crefformat{theo}{#2Theorem~\textup{#1}#3}
\crefformat{lem}{#2Lemma~\textup{#1}#3}
\crefformat{app}{#2Appendix~\textup{#1}#3}
\crefformat{scenario}{#2Scenario~\textup{#1}#3}
 \AtBeginDocument{%
    \crefname{figure}{Figure}{figures}%
}
\let\origref\cref
\def\cref#1{\textbf{\origref{#1}}}

\expandafter\def\expandafter\UrlBreaks\expandafter{\UrlBreaks%
  \do\a\do\b\do\c\do\d\do\e\do\f\do\g\do\h\do\i\do\j%
  \do\k\do\l\do\m\do\n\do\o\do\p\do\q\do\r\do\s\do\t%
  \do\u\do\v\do\w\do\x\do\y\do\z\do\A\do\B\do\C\do\D%
  \do\E\do\F\do\G\do\H\do\I\do\J\do\K\do\L\do\M\do\N%
  \do\O\do\P\do\Q\do\R\do\S\do\T\do\U\do\V\do\W\do\X%
  \do\Y\do\Z}

\usepackage[sort&compress,square,comma,numbers]{natbib}
\usepackage{fullpage,authblk}

\begin{document}

\title{A Time-Constrained Capacitated Vehicle Routing Problem in Urban E-Commerce Delivery}

\author[a,b]{Taner Cokyasar}
\author[a]{Anirudh Subramanyam}
\author[a]{Jeffrey Larson}
\author[a]{Monique Stinson}
\author[a]{Olcay Sahin}

\affil[a]{Argonne National Laboratory, 9700 S. Cass Ave. Lemont, IL, 60439, USA. {\tt tcokyasar@anl.gov}}
\affil[b]{Tarsus University, Takbas mah., Kartaltepe Sk. Tarsus, Mersin, 33400, Turkey}
\date{}                     %
\renewcommand\Affilfont{\itshape\small}

\maketitle

\begin{abstract}
Electric vehicle routing problems can be particularly complex when recharging must be performed mid-route.  
In some applications such as the e-commerce parcel delivery truck routing, however, mid-route recharging may not be necessary 
because of constraints on vehicle capacities 
and maximum allowed time for delivery. 
In this study, we develop a mixed-integer optimization model that exactly solves such a time-constrained capacitated vehicle routing problem, especially of interest to e-commerce parcel delivery vehicles. 
We compare our solution method with an existing metaheuristic and
carry out exhaustive case studies considering four U.S. cities---Austin, TX; Bloomington, IL; Chicago, IL; and Detroit, MI---and 
two vehicle types: conventional vehicles and battery electric vehicles (BEVs). 
In these studies we examine the impact of vehicle capacity, maximum allowed travel time, 
service time (dwelling time to physically deliver the parcel), and BEV range on system-level performance metrics including vehicle miles traveled (VMT).
We find that the service time followed by the vehicle capacity plays a key role in the performance of our approach. 
We assume an 80-mile BEV range as a baseline without mid-route recharging. Our results show that BEV range has a minimal impact on performance metrics because the VMT per vehicle averages around 72 miles. 
In a case study for shared-economy parcel deliveries, we observe that 
VMT could be reduced by 38.8\% in Austin if service providers were to operate their distribution centers jointly.

\hfill\break%
\noindent\textit{Keywords}: Time-constrained vehicle routing, Delivery planning, Optimization
\end{abstract}

\section{Introduction}
Vehicle routing problems (VRPs) are NP-hard problems that are fundamental in the 
transportation science field~\cite{archetti2011complexity}. Solving a VRP
requires determining optimal routes for a set of vehicles so that each 
location in a set of places is visited at least once. Naturally, many
VRP variants exist.
A time-constrained (vehicle-load) capacitated VRP (TCVRP) is an
important problem variant that is similar to the well-studied
distance-constrained VRP (DVRP)~\cite{almoustafa2013new,kara2007energy,kara2011arc,kek2008distance}.
The TCVRP considers optimally routing vehicles through a network 
to deliver packages to a set of locations subject to constraints on the
total travel time and the number of packages delivered by each vehicle.
In small VRP instances (e.g., tens of delivery locations and vehicles), optimal
solutions can be identified in a reasonable amount of time~\cite{arnold2019efficiently}. These
routing problems become challenging at large scales with hundreds of
thousands of delivery locations and multiple depots (the unique starting and
ending location for subsets of vehicles), although numerous heuristic and metaheuristic solution approaches exist in the literature. 
In this study, we formulate a mixed-integer program (MIP) to exactly solve small (e.g., 50 customers) TCVRP instances. 
Using validated simulation data for four cities, we conduct case studies investigating the impact of battery electric vehicles (BEVs) 
on energy consumption compared with conventional vehicles (CVs) in e-commerce parcel deliveries at an urban scale. 
We carry out sensitivity analyses to highlight the importance of service (i.e., package dropping) times and 
to determine whether BEV ranges play a role in the energy consumption of parcel delivery trucks.

Large-scale VRPs appear in many real-world and simulated transportation networks. 
Our  work here is motivated by a study of the effects of optimal
delivery truck tours in POLARIS, the Planning and Operations
Language for Agent-based Regional Integrated Simulation~\cite{auld2016polaris}. This software is
frequently used to
quantify the impact of emerging and existing vehicle and transportation
technologies on a variety of metrics, such as vehicle miles traveled (VMT),
energy consumed, and greenhouse gas emitted in large metropolitan areas. 
VRPs abound within POLARIS, but
a common instance that is increasingly important to model accurately is the
effect of package delivery (from Amazon, FedEx, UPS, USPS, etc.) at the system level.
Solving the truck routing problem at a
large scale allows estimating an average VMT per vehicle, which then informs what BEV range to be satisfactory in this application. Furthermore, the energy consumption of BEVs and CVs can be estimated to quantify the marginal benefit of using BEVs at a system level.

Compared with CV routing, BEV routing---namely, electric VRP (EVRP)---is complex because of the en-route charging need. 
Travel time to arrive at a charging station, waiting time due to congestion at a station, 
time to recharge, and when to recharge complicate the EVRP. 
Apart from these factors, the EVRP models are similar to VRP models. 
In this study, we consider a case where delivery BEVs leave a designated depot fully charged, 
make deliveries to customers, and return to the depot before  running out of battery. 
Under such a setting, BEVs are not allowed en-route charging, 
and hence the problem becomes a TCVRP in which only vehicle capacities and service times are constrained. 
To account for the BEV distance range constraint, we use methodologies developed in the DVRP literature.

The contribution of this study is quantifying energy consumption (as a linear function of the VMT) of e-commerce delivery BEVs and CVs at a regional scale for large metropolitan areas supported by validated simulation data under various conditions (e.g., BEV range, service time, vehicle capacity, and work hours). By carrying out sensitivity analyses, We consider improvement in these conditions and reveal the ones that need more attention to improve the system-level performance metrics, such as VMT, vehicle hours traveled (VHT), and the number of vehicles needed. For instance, our analyses show that the service time and vehicle capacity (maximum number of packages vehicles can deliver in a route) are the key determiners to improve the metrics. Moreover, we provide managerial insights into the cases in which the BEV range is an impactful factor. Although we provide an MIP model to solve the TCVRP, most of its components existed in the literature. Therefore, the key contribution is the application of this approach in an e-commerce delivery context to convey insights using validated simulation data for various areas.

\section{Literature Review}
Research in VRPs started in earnest with the 1959 paper ``The Truck
Dispatching Problem" of Dantzig and Ramser~\cite{dantzig1959truck}. The authors
introduced the problem in detail and highlighted its resemblance to
the traveling salesman problem (TSP) studied in~\cite{flood1956traveling}. Since then, numerous variants
of the problem have been studied, and alternative solution approaches have been 
proposed~\cite{bertsimas1992vehicle,ordonez2007priori,ehmke2018optimizing,figliozzi2010vehicle,drexl2012synchronization,wang2013vehicle}.
For further information, see recent surveys of the VRP 
literature~\cite{eksioglu2009vehicle,
braekers2016vehicle,konstantakopoulos2020vehicle}.

We study the TCVRP with asymmetric travel costs, that is, when the cost of traveling
from some location A to location B may not be the same as the cost of traveling
from B to A. This asymmetry is a result of unidirectional links in
the transportation network; the literature commonly uses the acronym ADVRP
(asymmetric distance-constrained vehicle routing problem) for versions of this
problem that do not consider vehicle
capacities~\cite{laporte1987branch,almoustafa2013new}. We use a flow-based ADVRP formulation
introduced in~\cite{kara2007energy} as an exact solution method by extending it to include vehicle 
capacity constraints. 
Although the resulting model can prove optimality for
a set of deliveries of a depot, doing so can require considerable 
computational resources. 
On the other hand, both the literature and various open-source platforms contain 
numerous heuristic and metaheuristic techniques to solve almost any type 
of VRPs to reasonable optimality bounds. 

While the asymmetric TCVRP would seem to be the most natural model for a
modern package delivery problem, relatively few studies can  be found in the
literature~\cite{almoustafa2013new,hashimoto2006vehicle,kara2007energy,kara2011polynomial}. Yet, this is not surprising because possible methodological improvement to the VRP is limited, and existing solution methods can be adjusted to account for various emerging aspects of the problem. In~\cite{almoustafa2013new}, the authors introduced an exact solution procedure for the ADVRP that can solve instances with 1,000 customers. A similar arc-based formulation to the one presented in this study was developed in~\cite{kara2011polynomial} to solve a distance- and capacity-constrained VRP. The difference in this study is that routes are time constrained. In~\cite{hashimoto2006vehicle}, the authors studied a VRP with flexible time-windows and travel times, proved that the problem is NP-hard, and used local search to solve the problem. One difference of their study is that they considered time-windows, that is each customer should be visited within a predefined time interval. In our problem, however, we consider a single strict time constraint for the last customer visited in a route, that is the last customer in a route should be served exactly at a certain time. In other words, the total time spent to serve all customers in a route should not exceed a predefined time that can be interpreted as work hours. Therefore, the VRP with time-windows (VRPTW) can also be viewed as a relevant literature. See~\cite{braysy2005vehicle,braysy2005vehicle2,kallehauge2008formulations} for further information about VRPTW.

Apart from the problem type considered,  we  also review the BEV routing literature. EVRPs are centered on en-route charging and battery swapping~\cite{chen2016electric, keskin2018matheuristic, loffler2020routing}. In~\cite{keskin2018matheuristic}, the authors considered an EVRP with time-windows (EVRPTW) and fast charging. They developed two mathematical models and tested them on small and large problem instances. Since typical delivery routes do not require more than one recharge, heuristic methods were developed to solve the EVRPTW on a single charge~\cite{loffler2020routing}. The  average FedEx VMT in the U.S. parcel deliveries was reported as 41.4 miles~\cite{barnitt2011fedex,feng2013economic}. Therefore, an en-route recharge may not be necessary for e-commerce deliveries in the real world. Our study simplifies the problem and assumes that vehicles are not recharged en-route and  that their routes are formed such that they can complete a route without the need for a recharge. Different scenarios comparing the routing of BEVs and CVs were studied in~\cite{davis2013methodology}. The authors analytically estimated the average cost of serving routes using a continuous approximation of the VRP rather than solving it. They concluded that high VMT, frequent stops at customers, and tax incentives make BEVs competitive in the long term. For a comprehensive review, see~\cite{kucukoglu2021electric}.

\section{Methodology}\label{sec:methodology}
We now describe the TCVRP in detail and present a solution
approach. Let the graph $\sG=\left(\sV, \sE\right)$
represent a network, where $\sV$ is the set of vertices and $\sE$ is the
set of arcs. 
Vertex $0$ denotes a depot from where vehicles are deployed and need to return at the end of a planning horizon, typically one day. Therefore, we use $\sV^\prime=\sV\setminus{\{0\}}$ to denote a set of customer locations. 
Let $Q$ and $\overline{T}$ be the capacity (maximum number of packages) and maximum allowed total 
travel time for each vehicle, respectively. Let $T_{ij}$ and $D_{ij}$
represent the travel time and the travel distance on arc $(i,j)\in\sE$, respectively. Let $S_i$ be the service time (also referred to as the dwell time~\cite{chen2014preliminary}) to be spent at vertex $i$. The parameter $N_i$ indicates the number of packages delivered at vertex $i\in\sV^\prime$. The binary 
variable $x_{ij}$ indicates 
whether an arc $(i,j)$ is 
traversed by a vehicle; if so,  
$x_{ij}=1$. We assume that
the number of vehicles is a variable denoted by $k$. To track 
the number of packages delivered at vertex $i$ while en route 
to $j$ (after leaving $i$ and $i\neq j$), we define $y_{ij}\in\mathbb{R}_{\geq 0}$. 
Similarly, we define $z_{ij}\in\mathbb{R}_{\geq 0}$ to track the total travel time from the depot to vertex $j$, where $i$ is the predecessor of $j$. 

The TCVRP is to route delivery vehicles so 
that their total travel distance is minimized while satisfying travel time and
vehicle capacity constraints on each vehicle. The travel time includes both the time spent traveling on arcs and the service time that is needed to park a vehicle and physically conduct a delivery. At this point
we may formulate an MIP to solve the TCVRP that minimizes the total travel time. Sets, parameters, and variables used in this
section are provided in \cref{DL_TCVRP_notations}. We attempt
to follow the notation used in the model of~\cite{almoustafa2013new} that we
are extending to include vehicle capacity constraints. The MIP to solve the TCVRP is as follows:

\begin{table}[!htb]
  \footnotesize
  \caption{Sets, parameters, and variables used in the depot-level
  TCVRP.}\label[tab]{DL_TCVRP_notations}
  \begin{tabularx}{\linewidth}{lX}
  \toprule
    \textbf{Set} & \textbf{Definition}\\
    \midrule
    $\sE$ & a set of arcs that can be traversed, indexed by $(i,j)$.\\
    $\sV$ & a set of vertices including a depot and customer locations, indexed by $0$, $i$, or $j$.\\
    $\sV^\prime$ & a subset of vertices representing customer locations to be visited; $\{0\}\cup\sV^\prime = \sV$, where vertex $0$ denotes the depot location.\\
    \midrule
    \textbf{Param.} & \textbf{Definition}\\
    \midrule
    $\overline{D}$ & maximum allowed travel distance for each vehicle.\\
    $D_{ij}$ & travel distance on arc $(i,j)$.\\
    $N_i$ & number of packages to be delivered at vertex $i$.\\
    $Q$ & package capacity of a vehicle.\\
    $S_i$ & service time (i.e., dwell time) at vertex $i$.\\
    $\overline{T}$ & maximum allowed travel time for each vehicle.\\
    $T_{ij}$ & travel time on arc $(i,j)$.\\
    \midrule
    \textbf{Var.} & \textbf{Definition}\\
    \midrule
    $k$ & number of vehicles to be used.\\
    $x_{ij}$ & $\begin{cases}
          1, & \text{if a vehicle drives on arc $(i,j)\in\sE,~i\neq j$},\\
          0, & \text{otherwise}.\\
         \end{cases}$\\
    $y_{ij}$ & number of packages delivered at vertex $i$ while en route to vertex $j$, i.e., after leaving $i$, where $i\neq j$.\\
    $z_{ij}$ & total travel time from the depot to vertex $j$, where $i$ is the predecessor of $j$ and $i\neq j$.\\ 
    $z_{ij}^\prime$ & total travel distance from the depot to vertex $j$, where $i$ is the predecessor of $j$ and $i\neq j$.\\ 
    \bottomrule
  \end{tabularx}
\end{table}

\begin{equation}\label[objfunc]{DL_TCVRP_obj_fun}
    \min_{\textbf{k},\textbf{x},\textbf{y},\textbf{z}} \sum_{(i,j)\in \sE}D_{ij}x_{ij},
\end{equation}
\noindent subject to,
\begin{equation}\label[consset]{DL_TCVRP_cons1}
    \sum_{i\in\sV} x_{ij} = 1\qquad \forall j\in \sV^\prime,
\end{equation}
\begin{equation}\label[consset]{DL_TCVRP_cons2}
    \sum_{j\in\sV} x_{ij} = 1\qquad \forall i\in \sV^\prime,
\end{equation}
\begin{equation}\label[cons]{DL_TCVRP_cons3}
    \sum_{i\in\sV^\prime} x_{0i} = k,
\end{equation}
\begin{equation}\label[cons]{DL_TCVRP_cons4}
    \sum_{i\in\sV^\prime} x_{i0} = k,
\end{equation}
\begin{equation}\label[consset]{DL_TCVRP_cons10}
    y_{ij} \leq Q x_{ij} \qquad \forall (i,j)\in \sE,
\end{equation}
\begin{equation}\label[consset]{DL_TCVRP_cons11}
    \sum_{j\in\sV} y_{ij}-\sum_{j \in\sV}y_{ji} = N_i \quad \forall i\in \sV^\prime,
\end{equation}
\begin{equation}\label[consset]{DL_TCVRP_cons5}
    \sum_{\substack{j\in \sV}} z_{ij} - \sum_{\substack{j\in \sV}}z_{ji} = \sum_{\substack{j\in \sV}} \left(T_{ij} + S_i \right)x_{ij} \qquad \forall i\in \sV^\prime,
\end{equation}
\begin{equation}\label[consset]{DL_TCVRP_cons6}
    z_{ij} \leq \left(\overline{T} - T_{j0}\right)x_{ij} \qquad \forall i\in V,j\in\sV^\prime,
\end{equation}
\begin{equation}\label[consset]{DL_TCVRP_cons7}
    z_{ij} \geq \left(T_{ij} + T_{0i} + S_i\right)x_{ij}  
    \qquad \forall i\in\sV^\prime,j\in\sV,
\end{equation}
\begin{equation}\label[consset]{DL_TCVRP_cons8}
    z_{i0} \leq \overline{T} x_{i0} \qquad \forall i\in \sV^\prime,
\end{equation}
\begin{equation}\label[consset]{DL_TCVRP_cons9}
    z_{0i} = T_{0i} x_{0i} \qquad \forall i\in \sV^\prime,
\end{equation}
\begin{equation*}
    \nonumber k\in\mathbb{Z}_{\geq 0},~x_{ij}\in\{0,1\},~ y_{ij},z_{ij}\in \mathbb{R}_{\geq 0}\qquad\forall (i,j)\in\sE.
\end{equation*}

Objective function \labelcref{DL_TCVRP_obj_fun} minimizes the total travel distance on arcs. Constraints \labelcref{DL_TCVRP_cons1}--\labelcref{DL_TCVRP_cons4} satisfy
the connectivity of the vehicle routes and are standard VRP constraints.
Constraints \labelcref{DL_TCVRP_cons10}--\labelcref{DL_TCVRP_cons11}
impose the vehicle capacity limitations, that is the number of packages delivered by each vehicle does not exceed $Q$. Constraints \labelcref{DL_TCVRP_cons5}--\labelcref{DL_TCVRP_cons9} ensure that the total travel time for
each vehicle does not exceed $\overline{T}$. (These constraints are
illustrated and explained in~\cite{kara2011arc}.)

Although total travel time is a natural constraint for CVs due to limited work hours, we need to further impose total distance constraints to consider BEV range limitations. Following the above model structure, this process is  straightforward. Let $\overline{D}$ represent the maximum allowed travel distance for each vehicle, and let $z_{ij}^\prime\in\mathbb{R}_{\geq 0}$ denote the total travel distance from the depot to vertex $j$, where $i$ is the predecessor of $j$ satisfying $i\neq j$. We may additionally introduce Constraints \labelcref{DL_TCVRP_cons12}--\labelcref{DL_TCVRP_cons16} to account for BEV range limitations.
\begin{equation}\label[consset]{DL_TCVRP_cons12}
    \sum_{\substack{j\in \sV}} z_{ij}^\prime - \sum_{\substack{j\in \sV}}z_{ji}^\prime = \sum_{\substack{j\in \sV}} D_{ij}x_{ij} \qquad \forall i\in \sV^\prime,
\end{equation}
\begin{equation}\label[consset]{DL_TCVRP_cons13}
    z_{ij}^\prime \leq \left(\overline{D} - D_{j0}\right)x_{ij} \qquad \forall i\in V,j\in\sV^\prime,
\end{equation}
\begin{equation}\label[consset]{DL_TCVRP_cons14}
    z_{ij}^\prime \geq \left(D_{ij} + D_{0i}\right)x_{ij}  
    \qquad \forall i\in\sV^\prime,j\in\sV,
\end{equation}
\begin{equation}\label[consset]{DL_TCVRP_cons15}
    z_{i0}^\prime \leq \overline{D} x_{i0} \qquad \forall i\in \sV^\prime,
\end{equation}
\begin{equation}\label[consset]{DL_TCVRP_cons16}
    z_{0i}^\prime = D_{0i} x_{0i} \qquad \forall i\in \sV^\prime
\end{equation}

Constraints \labelcref{DL_TCVRP_cons12}--\labelcref{DL_TCVRP_cons16} function similarly to constraints \labelcref{DL_TCVRP_cons5}--\labelcref{DL_TCVRP_cons9}. We note that the majority of the model components have already existed in the literature;  our  contribution is the addition of Constraints \labelcref{DL_TCVRP_cons10}--\labelcref{DL_TCVRP_cons11} to the MIP presented in~\cite{kara2011arc}.

\section{Case Studies}\label{sec:num_exp} 
In this section, we first thoroughly explain the design of experiments, laying out all implementation details. We then describe computational experiments that show the quality of the solution method, and we compare it with an iterated tabu search (ITS) metaheuristic from the literature. We provide extensive sensitivity analyses to compare BEVs and CVs under various cases focusing on three large U.S. cities---Austin, TX; Chicago, IL; and Detroit, MI---and a small city, Bloomington, IL.

\subsection{Design of Experiments}\label{sec:experiments}

The POLARIS agent-based modeling framework was used to generate problem
instances for four cities: Austin, Bloomington, 
Chicago, and Detroit~\cite{auld2016polaris}. (We use the words \textit{city} and \textit{area} interchangeably except when we refer to the physical area of a region.) 
The National Household Travel Survey revealed that a household places approximately one order per week~\cite{nhts}. Hence, POLARIS assumes that nearly 1/7 of households (which we also refer to as customers) require an e-commerce delivery service on a typical day. It randomly draws their locations from the databases following a uniform distribution, that is each household has an equal probability to be selected. Since POLARIS assumes some customers do not receive a delivery service, they are exempted even if selected randomly. \cref{problem_size} tabulates network topology and other parameters for these cities. The first column denotes the number of households in the area. The second column indicates the number of households to be delivered to on the considered day. (We assume each household requests one delivery, although the assumption can be easily relaxed.) The third and the fourth columns show the number of arcs and vertices in the area's road network, respectively. Here, arcs refer to unidirectional road segments in the road network, and vertices are connectors that are on both ends of arcs. Arcs and vertices are used to compute the shortest paths between any given points. The fifth and the sixth columns show the number of e-commerce delivery centers (i.e., depots) and service providers servicing in these areas. These providers are Amazon, FedEx, UPS, and USPS. We identified the number of depots for these providers in the four areas from publicly available sources. Depots of these providers have different uses. Although we paid attention to include only the depots that allow parcel deliveries, some depots may not really serve this purpose in reality. In Bloomington, we could not locate any Amazon depots and hence considered the three providers. 

\begin{table*}[!htbp]
\footnotesize
\centering
\caption{Network topology and other parameters
for the areas considered in experiments.}\label[tab]{problem_size}
{\begin{tabular*}{\linewidth}{l@{\extracolsep{\fill}}ccccccc}
\toprule
City & Area (sq mi) & \# households & \# households ordering & \# arcs & \# vertices & \# depots & \# providers\\
\midrule
Austin & 5,377 & 830,000 & 158,172 & 40,891 & 17,231 & 22 & 4 \\
Bloomington & 74 & 16,605 & 2,816 & 7,013 & 2,540 & 8 & 3\\
Chicago & 11,116 & 4,017,583 & 606,669 & 57,267 & 19,377 & 53 & 4 \\
Detroit & 4,635 & 1,910,260 & 271,129 & 60,701 & 26,424 & 30 & 4\\
\bottomrule
\end{tabular*}}
\end{table*}

\begin{figure*}[!htb]
\centering
\subfloat[\centering Austin, TX.\label{fig:Austin_map}]{%
\includegraphics*[width=0.48\textwidth,height=\textheight,keepaspectratio]{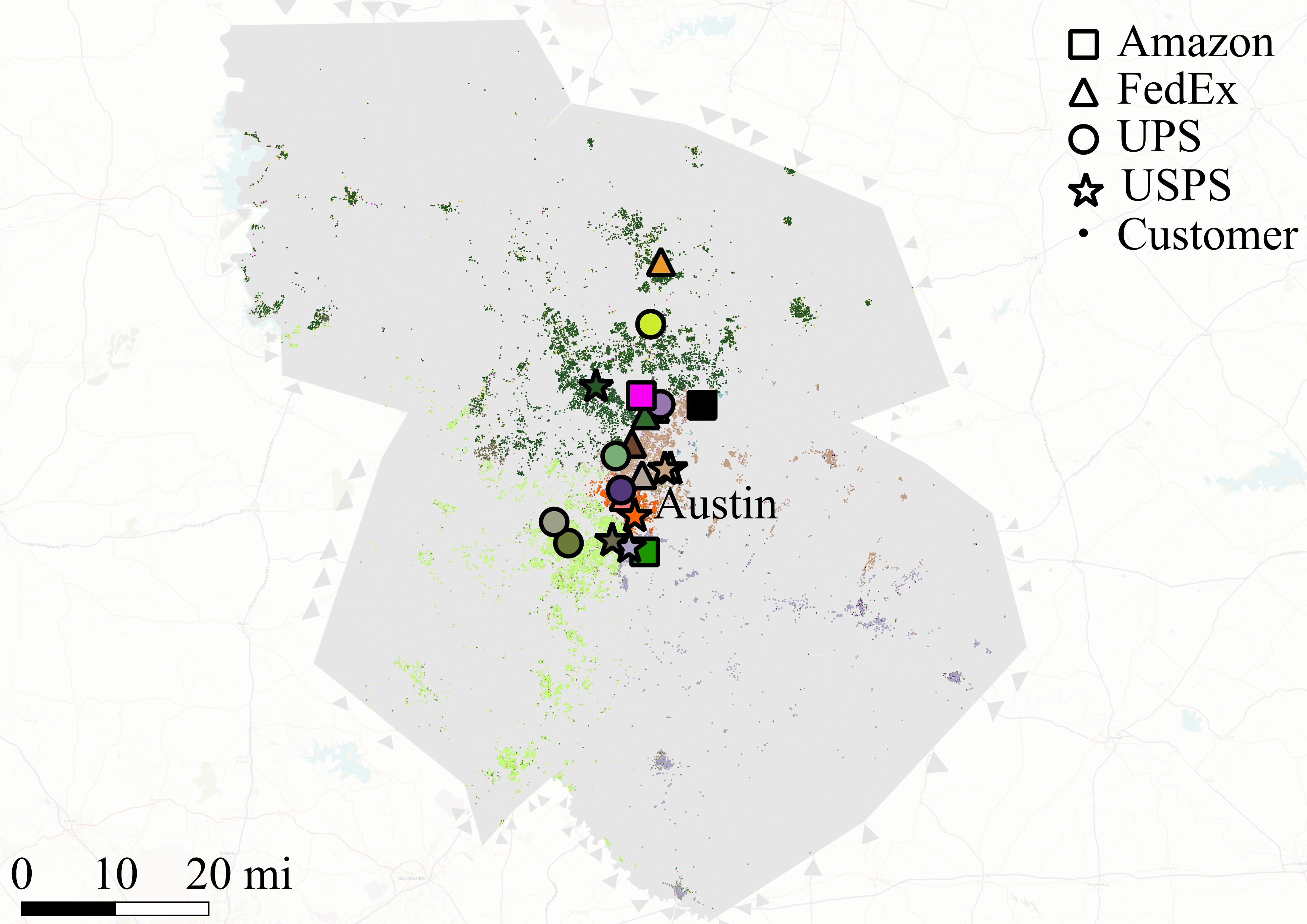}}
\quad
\subfloat[\centering Bloomington, IL.\label{fig:Bloomington_map}]{%
\includegraphics*[width=0.48\textwidth,height=\textheight,keepaspectratio]{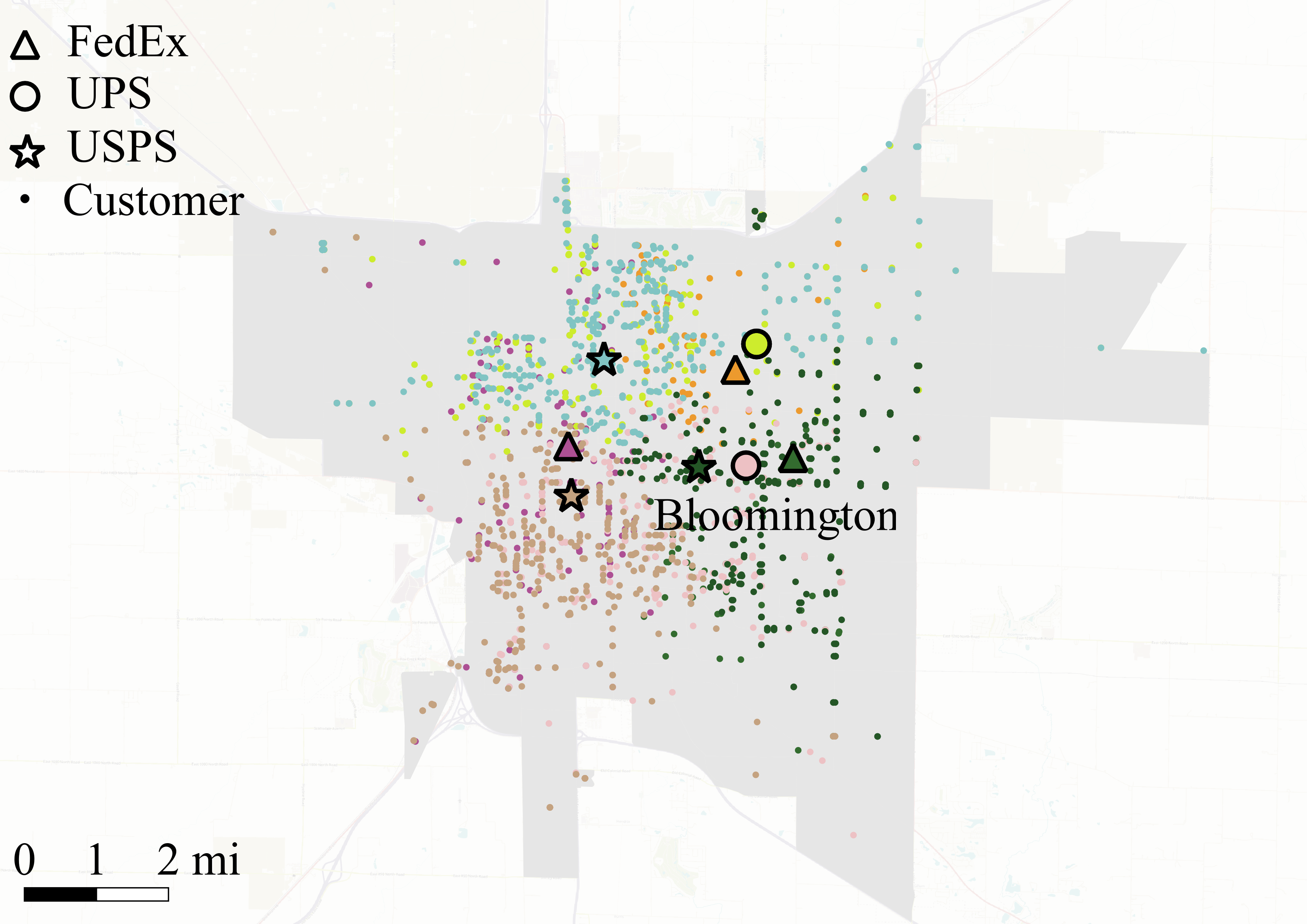}}\\
\subfloat[\centering Chicago, IL.\label{fig:Chicago_map}]{%
\includegraphics*[width=0.48\textwidth,height=\textheight,keepaspectratio]{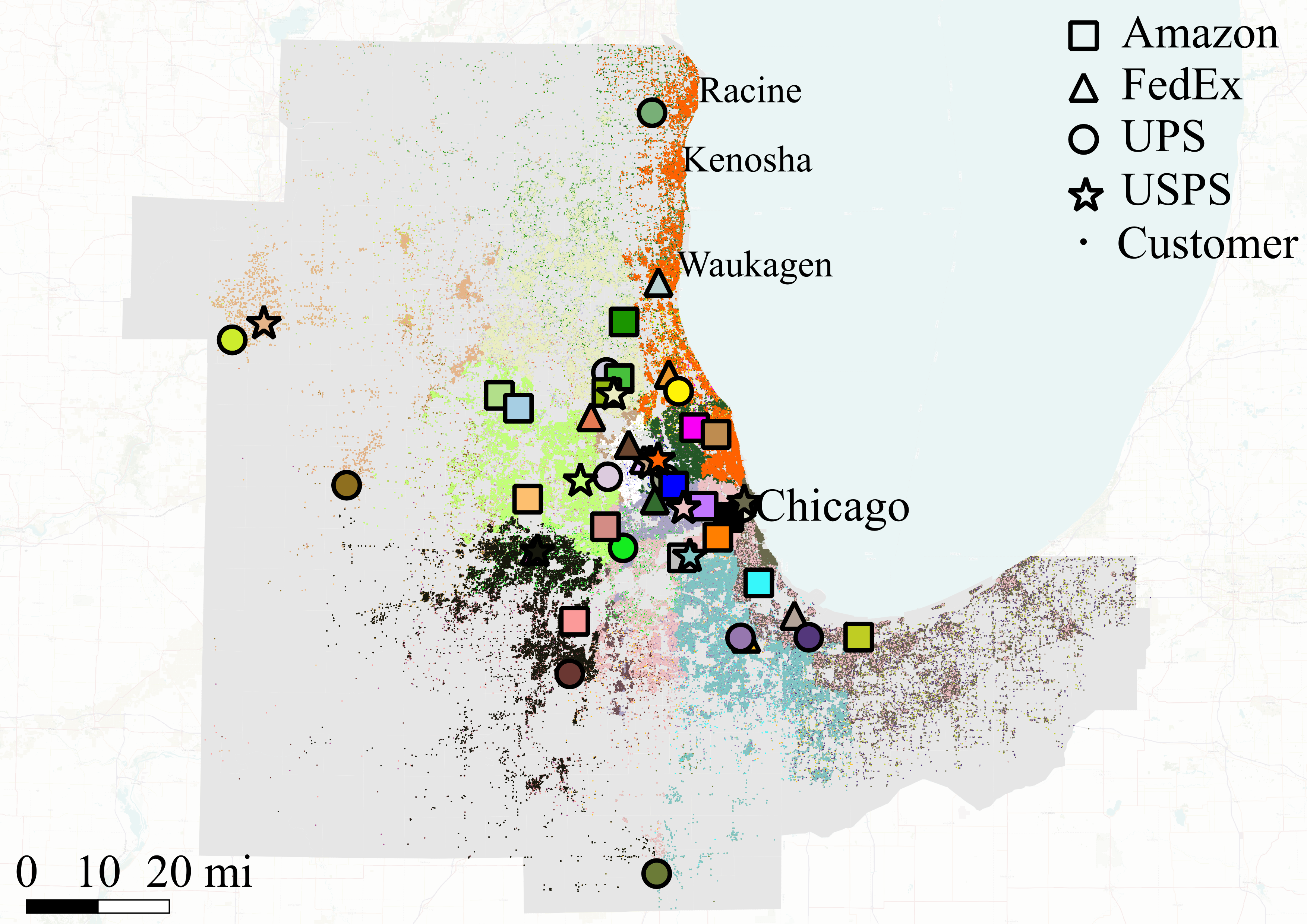}}
\quad
\subfloat[\centering Detroit, MI.\label{fig:Detroit_map}]{%
\includegraphics*[width=0.48\textwidth,height=\textheight,keepaspectratio]{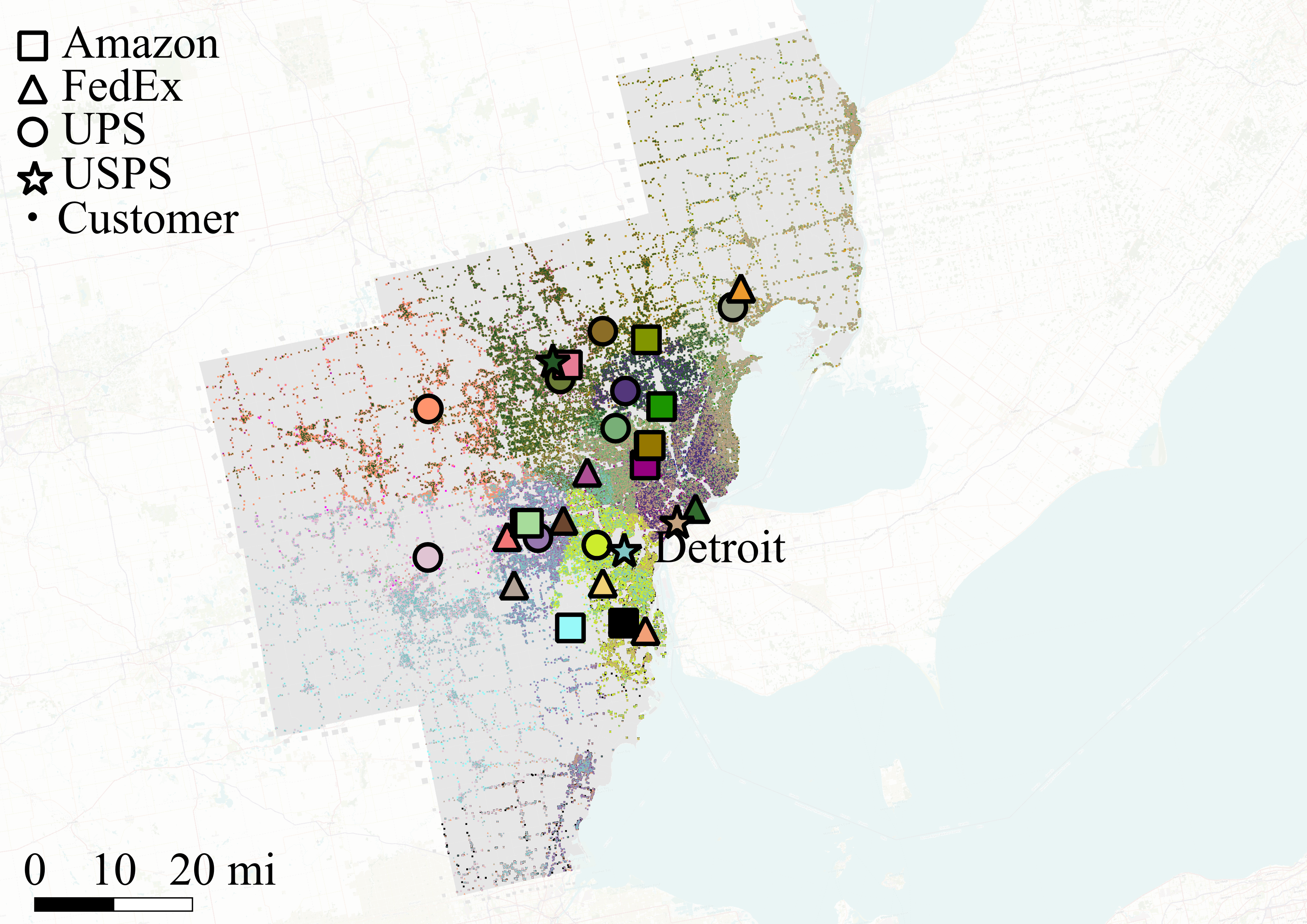}}
\caption{Illustrative problem layouts. Customer locations are color-coded to match the colors of their assigned depots. Certain colors dominate the maps because of overlapping points.} \label[fig]{maps}
\end{figure*}

In the three cities, we randomly distributed customers to Amazon, FedEx, UPS, and USPS following 21, 16, 24, and 39 percentage shares, respectively~\cite{pitney}. 
In Bloomington, we equally distributed 21\% of Amazon's shares to the three providers.
Since solving the depot-to-customer assignments and the VRP in conjunction complicates the problem, 
we assume that providers solve assignment problems before the VRP to determine a set of customers to be served by each depot. 
The assignment problem minimizes the total travel distance between customers and depots while 
adhering to capacity and assignment constraints. The capacity here refers to a limit on the number of 
customers to be assigned to a single depot, and the assignment constraint ensures that each customer is assigned to exactly one depot. 
We do not present this model because we consider it to be out of the scope for this study. See \cite{giosa2002new} for reference and \cite{drexl2015survey,laporte1986exact} for the location routing problem that incorporates both routing and depot location decisions in a single model. 
\cref{stats1} shows the statistics on the number of customers at depot-level problems in the 
four cities. For instance, Chicago depots have an average of 11,447 customers. 
Each depot-level problem is an instance of the TCVRP. \cref{maps} illustrates resulting problem layouts for each city and demonstrates the area with a gray background.

\begin{table}[!htb]
\footnotesize
\centering
\caption{Statistics on the number of customers at depot-level problems.}\label[tab]{stats1}
{\begin{tabular*}{\linewidth}{l@{\extracolsep{\fill}}cccc}
\toprule
City & Avg. & Min. & Max. & Std. dev.\\
\midrule
Austin & 7,190 & 242 & 24,000 & 5,950 \\
Bloomington & 352 & 167 & 480 & 116 \\
Chicago & 11,447 & 905 & 25,200 & 7,466 \\
Detroit & 9,037 & 2,138 & 14,400 & 2,144\\
\bottomrule
\end{tabular*}}
\end{table}

\cref{stats1} shows that some TCVRP instances are  large, for example, 25,200 customers in one of the Chicago TCVRP instances. 
We cannot solve such large problems using the proposed MIP: It is prohibitive to compute and store the travel time and distance matrices. 
More important, since POLARIS is a mesoscopic traffic simulation tool, it does not contain microscopic network details, such as street-level minor roads.
Instead, road networks within the tool are composed of interstates, principal and other arterials, and major collectors.
To simplify the problem, we aggregated customers at midpoints of arcs that we call super-locations. 
In a depot-level problem, we found the closest super-location to each customer and the depot. 
The set of these super-locations is equivalent to $\sV$. 
We computed shortest paths in terms of travel time using Dijkstra's algorithm with the network information in POLARIS (i.e., vertices, arcs, arc speeds, and arc lengths). 
These paths yielded $T_{ij}$ and $D_{ij}$ parameter values between all super-locations in a TCVRP instance. 

\begin{figure}[!htb]
    \begin{center}
        \includegraphics[width=\linewidth]{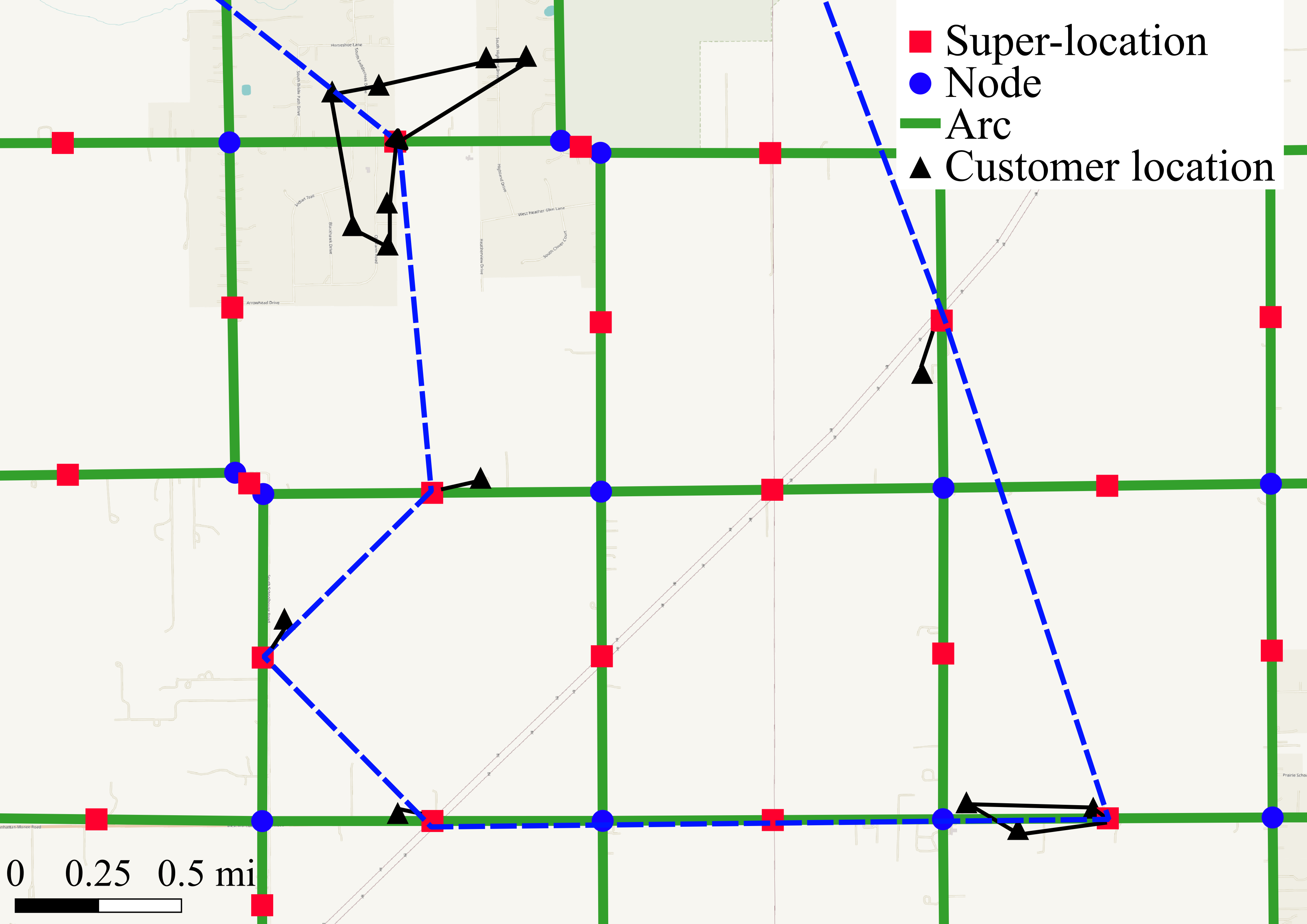}
    \end{center}
    \caption{An example network showing the relation of super-locations with road network components.\label[fig]{super_location}}
\end{figure}

We assume each customer requests one package; therefore, the total number of customers to be served at each super-location is equivalent to $N_i$ in the MIP. 
Since minor road data were unavailable and trucks drive at a low speed on minor roads, 
we estimated the customer-to-customer travel times by dividing Manhattan distances between customers by a constant speed of 15 mph. 
We solved TSPs with an objective of travel time minimization to optimize the sequence of visits at each super-location. 
We used the local search and the simulated annealing metaheuristics of the open-source python-tsp library~\cite{pythontsp}.
Best solutions obtained from these approximations were then pushed into Gurobi's TSP solver (modified to account for asymmetry) as a warm start~\cite{gurobi_tsp}.
Eventually, all TSPs were solved to optimality.
The value of $S_i$ in the MIP is the sum of the super-location-level travel time and a predetermined $P$ minutes per customer that accounts for the dwell time.
We also included the distance traveled at each super-location in the distance matrix entries, $D_{ij}$. \cref{super_location} illustrates an example road network layout along with super-locations, customer locations, and some routes. In this figure, black lines represent a super-location-level customer visit sequence and the blue dashed line shows the sequence of super-locations. Note that these lines do not show the real link driving patterns.

\cref{stats2} shows the statistics on the number of super-locations at depot level for each area. 
Since some instances are still  large, we use an ITS metaheuristic that was developed in a prior study. 
We refer the interested reader to~\cite{Subramanyam2020} for the 
implementation details of the ITS. To justify the quality of the ITS, 
we compare it with the MIP in the following section.

\begin{table}[!htbp]
\footnotesize
\centering
\caption{Statistics on the number of super-locations at depot-level problems.}\label[tab]{stats2}
{\begin{tabular*}{\linewidth}{l@{\extracolsep{\fill}}cccc}
\toprule
City & Avg. & Min. & Max. & Std. dev.\\
\midrule
Austin & 975 & 25 & 2,663 & 712 \\
Bloomington & 191 & 83 & 269 & 68 \\
Chicago & 1,346 & 93 & 3,707 & 872 \\
Detroit & 1,733 & 332 & 4,290 & 835\\
\bottomrule
\end{tabular*}}
\end{table}

Unless otherwise noted, we used $Q=120$ customers (i.e., \textit{packages} since every customer is assumed to request one package), $\overline{T}=10$ hours, and $P=2$ minutes for both BEVs and CVs; 
and we set $\overline{D}=80$ miles as the BEV range, as in~\cite{krok}. 
To account for the energy consumption of diesel-powered CVs in electricity units, 
we assume a CV fuel efficiency of eight mpg~\cite{eia}, where  one gallon of diesel is equivalent to 40.15 kWh energy~\cite{eia2}. 
We assume BEVs consume 1.14 kWh per mile~\cite{cerc}.
Although we do not calculate the energy consumption for all instances, 
these metrics can be used to observe the magnitude by multiplying the VMT by the kWh per mile.

All MIP computations were carried
out on an Intel{\textsuperscript \textregistered} Xeon{\textsuperscript \textregistered} Gold 6138
CPU@2.0 GHz workstation with 128 GB of RAM and 40 cores. Problem instances were solved by using
the Python 3.8.8 interface to the commercial solver Gurobi 9.1~\cite{gurobi}. 
For instances considered in the next section, we used Windows Subsystems for Linux to run the ITS (once per instance) on this workstation. 
For all other instances, throughout we used four workstations identical to the aforementioned, 
ran the ITS 10 times each with a one-hour time limit, and reported the best outcome of the 10 runs.

\subsection{Computational Performance of the MIP and the ITS methods}

We analyzed the computational performance of the MIP and compared it with the ITS method focusing solely on CVs. 
We designed small TCVRP instances by randomly sampling super-locations from the original four city problems that are based on POLARIS outputs. 
Let $V$ represent the number of super-locations.
From the datasets of each depot, we randomly drew $V\in\{25,~50,~100\}$. 
For each of these instances, we considered three scenarios. 
The first scenario assumes a baseline of $Q=60$ and $\overline{T}=8$ hours. 
In the second and the third, we set $Q=80$ and $\overline{T}=10$ hours, respectively. 
Since some depot-level problems have fewer than 100 super-locations,
we have slightly less than (number of depots) $\times$ 3 $\times$ 3 TCVRP instances for each city in total. 
We have made 24 of these problem instances (six instances for each city) and a
formulation of our MIP available at \url{https://gitlab.com/tcokyasar/tcvrp}. 
We carried out these analyses on all four cities to observe 
whether the outcomes were alike on different network configurations.
We capped the computational time at 300 seconds for instances with 
100 super-locations and 60 seconds for others in both the MIP solver and the ITS. 

\cref{performance} summarizes the performance of both methods. 
The second column aggregates the results on a scenario basis. 
It first reports the three scenario statistics separately, then lists average results for each $V$, and  shows the average results for all instances.
The third column denotes the number of instances in each scenario. 
Columns 4--7 and columns 8--10 categorize the results based on optimality and nonzero gap solutions. 
The MIP solver produces a lower bound, $\iota$, and an upper bound, $\upsilon$, 
for the objective value. The percent MIP gap is defined by $(1-\iota/\upsilon)\times 100$. 
The optimality condition is met when the MIP gap is below a default threshold of Gurobi. 
In the ITS instances, we calculate a percent ITS gap by comparing the best-found solution, $\omega$, 
with the MIP's $\iota$, that is, $(1-\omega/\iota)\times 100$. 
For this reason, the number of optimal instances of the ITS shown in column 5 
can be greater than the number of optimal instances of the MIP in column 4. 
Moreover, a positive ITS gap does not mean that the solution found by the method is not optimal 
because $\iota$ of the MIP is not guaranteed to be optimal.
Columns 6 and 7 indicate the average time to achieve optimality with the MIP 
solver and the ITS, respectively. In the nonzero gap portion of the table, we report the percent MIP gap followed by the percent ITS gap and the ITS time. 
Although the ITS runs during the whole allotted amount of time, 
the reported averages are based on the time when the best solutions are found. 
In the MIP case, however, the solver terminates once the optimality threshold is satisfied; it keeps running until the time limit is reached otherwise.

\begin{table*}[!htb]
    \footnotesize
    \centering
    \caption{Summary of computational performance of the MIP and the ITS.}\label[tab]{performance}
    {\begin{tabular*}{\linewidth}{c@{\extracolsep{\fill}}ccccccccc}
    \toprule
    \multicolumn{1}{c}{\multirow{2}{*}{City}} & \multicolumn{1}{c}{\multirow{2}{*}{Scenario}} & \multicolumn{1}{c}{\multirow{2}{*}{\# inst.}} & \multicolumn{4}{c}{Optimal} & \multicolumn{3}{c}{Nonzero Gap}\\ \cmidrule(l{0.5em}r{0.5em}){4-7} \cmidrule(l{0.5em}r{0.5em}){8-10}
    \multicolumn{1}{c}{}& \multicolumn{1}{c}{}& \multicolumn{1}{c}{}& \multicolumn{1}{c}{\# MIP inst.} & \multicolumn{1}{c}{\# ITS inst.} & \multicolumn{1}{c}{MIP time (s)} & \multicolumn{1}{c}{ITS time (s)} & \multicolumn{1}{c}{MIP gap (\%)} & \multicolumn{1}{c}{ITS gap (\%)} & \multicolumn{1}{c}{ITS time (s)} \\
    \midrule
    \parbox[t]{2mm}{\multirow{8}{*}{\rotatebox[origin=c]{90}{Austin}}} & 1 & 62 & 17 & 17 & 6.32 & 0.22 & 8.64 & 0.06 & 97.5 \\
    & 2 & 62 & 20 & 20 & 6.83 & 0.03 & 7.38 & 0.05 & 60.7 \\
    & 3 & 62 & 19 & 19 & 7.52 & 0.04 & 7.43 & 0.05 & 74.9 \\
    \cmidrule(l{0.5em}r{0.5em}){2-2}
    \cmidrule(l{0.5em}r{0.5em}){3-3}
    \cmidrule(l{0.5em}r{0.5em}){4-7} \cmidrule(l{0.5em}r{0.5em}){8-10}
    & $V=25$ & 66 & 55 & 55 & 6.6 & 0.05 & 4.4 & 0.04 & 0.04 \\
    & $V=50$ & 60 & 1 & 1 & 23.87 & 2.25 & 5.51 & 0.04 & 10 \\
    & $V=100$ & 60 & 0 & 0 & N/A & N/A & 10.75 & 0.07 & 159.4 \\
    \cmidrule(l{0.5em}r{0.5em}){2-2}
    \cmidrule(l{0.5em}r{0.5em}){3-3}
    \cmidrule(l{0.5em}r{0.5em}){4-7} \cmidrule(l{0.5em}r{0.5em}){8-10}
    & All & 186 & 56 & 56 & 6.91 & 0.09 & 7.83 & 0.05 & 78.1 \\
    \midrule
    \parbox[t]{2mm}{\multirow{8}{*}{\rotatebox[origin=c]{90}{Bloomington}}} & 1 & 23 & 17 & 17 & 26.9 & 15.4 & 4.1 & 0.03 & 112 \\
    & 2 & 23 & 17 & 19 & 39.1 & 51.8 & 1.4 & 0.01 & 40 \\
    & 3 & 23 & 18 & 17 & 19.6 & 0.8 & 2.19 & 0.01 & 169 \\
    \cmidrule(l{0.5em}r{0.5em}){2-2}
    \cmidrule(l{0.5em}r{0.5em}){3-3}
    \cmidrule(l{0.5em}r{0.5em}){4-7} \cmidrule(l{0.5em}r{0.5em}){8-10}
    & $V=25$ & 24 & 23 & 24 & 0.9 & 0.01 & 0.01 & N/A & N/A \\
    & $V=50$ & 24 & 20 & 21 & 14.5 & 1.5 & 1.65 & 0.02 & 3.1 \\
    & $V=100$ & 21 & 9 & 8 & 129 & 154 & 3.11 & 0.02 & 141 \\
    \cmidrule(l{0.5em}r{0.5em}){2-2}
    \cmidrule(l{0.5em}r{0.5em}){3-3}
    \cmidrule(l{0.5em}r{0.5em}){4-7} \cmidrule(l{0.5em}r{0.5em}){8-10}
    & All & 69 & 52 & 53 & 28.3 & 23.8 & 2.59 & 0.02 & 115 \\
    \midrule
    \parbox[t]{2mm}{\multirow{8}{*}{\rotatebox[origin=c]{90}{Chicago}}} & 1 & 158 & 63 & 65 & 12.1 & 0.05 & 8.31 & 0.05 & 69.9 \\
    & 2 & 158 & 66 & 68 & 11.2 & 0.04 & 7.84 & 0.05 & 86.0 \\
    & 3 & 158 & 68 & 70 & 12.4 & 0.07 & 7.71 & 0.05 & 65.7 \\
    \cmidrule(l{0.5em}r{0.5em}){2-2}
    \cmidrule(l{0.5em}r{0.5em}){3-3}
    \cmidrule(l{0.5em}r{0.5em}){4-7} \cmidrule(l{0.5em}r{0.5em}){8-10}
    & $V=25$ & 159 & 151 & 157 & 10.3 & 0.05 & 2.25 & 0.03 & 0.31 \\
    & $V=50$ & 159 & 46 & 46 & 17.2 & 0.08 & 6.47 & 0.04 & 12.9 \\
    & $V=100$ & 156 & 0 & 0 & N/A & N/A & 9.31 & 0.06 & 119.1 \\
    \cmidrule(l{0.5em}r{0.5em}){2-2}
    \cmidrule(l{0.5em}r{0.5em}){3-3}
    \cmidrule(l{0.5em}r{0.5em}){4-7} \cmidrule(l{0.5em}r{0.5em}){8-10}
    & All & 474 & 197 & 203 & 11.9 & 0.05 & 7.95 & 0.05 & 73.9 \\
    \midrule
    \parbox[t]{2mm}{\multirow{8}{*}{\rotatebox[origin=c]{90}{Detroit}}} & 1 & 90 & 30 & 31 & 8.08 & 0.78 & 9.41 & 0.06 & 79.0 \\
    & 2 & 90 & 34 & 34 & 9.48 & 0.58 & 8.57 & 0.05 & 86.5 \\
    & 3 & 90 & 26 & 26 & 4.84 & 0.01 & 9.05 & 0.06 & 81.7 \\
    \cmidrule(l{0.5em}r{0.5em}){2-2}
    \cmidrule(l{0.5em}r{0.5em}){3-3}
    \cmidrule(l{0.5em}r{0.5em}){4-7} \cmidrule(l{0.5em}r{0.5em}){8-10}
    & $V=25$ & 90 & 84 & 84 & 5.4 & 0.02 & 3.58 & 0.04 & 0.01 \\
    & $V=50$ & 90 & 6 & 7 & 39.9 & 6.2 & 6.09 & 0.04 & 8.4\\
    & $V=100$ & 90 & 0 & 0 & N/A & N/A & 12.12 & 0.07 & 156 \\
    \cmidrule(l{0.5em}r{0.5em}){2-2}
    \cmidrule(l{0.5em}r{0.5em}){3-3}
    \cmidrule(l{0.5em}r{0.5em}){4-7} \cmidrule(l{0.5em}r{0.5em}){8-10}
    & All & 270 & 90 & 91 & 7.7 & 0.5 & 9.02 & 0.06 & 82.3 \\
    \bottomrule
    \multicolumn{10}{l}{\scriptsize{\textit{Note}: N/A = not applicable}} \\
    \end{tabular*}}
\end{table*}

\cref{performance} shows that the ITS method outperforms the MIP method in all scenarios in terms of the solution time. 
Both methods were effective on approximately the same number of instances (see the number of optimal instances for both methods). 
Only one scenario (Bloomington's $V=100$) provided a higher number of optimal solutions in the MIP compared with the ITS.
We can conclude that the MIP is not scalable to large problem instances 
and that the ITS provides better solutions than  the MIP does for most instances. 
On the other hand, the MIP provides a lower bound solution using which we can ensure a confidence interval for the solutions gained from the ITS. 
Therefore, an exact solution method---although not guaranteed to perform well at problem scales desired to be solved---is important 
to have on hand to assess the quality of alternative methods.
We employ the ITS method for all experiments henceforth.

\subsection{Impact of Vehicle Capacity on System-Level Metrics}
\begin{figure*}[!htbp]
    \centering
    \subfloat[\centering Austin\label{fig:Aus_givenQ}]{%
        \includegraphics*[width=0.25\textwidth,height=\textheight,keepaspectratio]{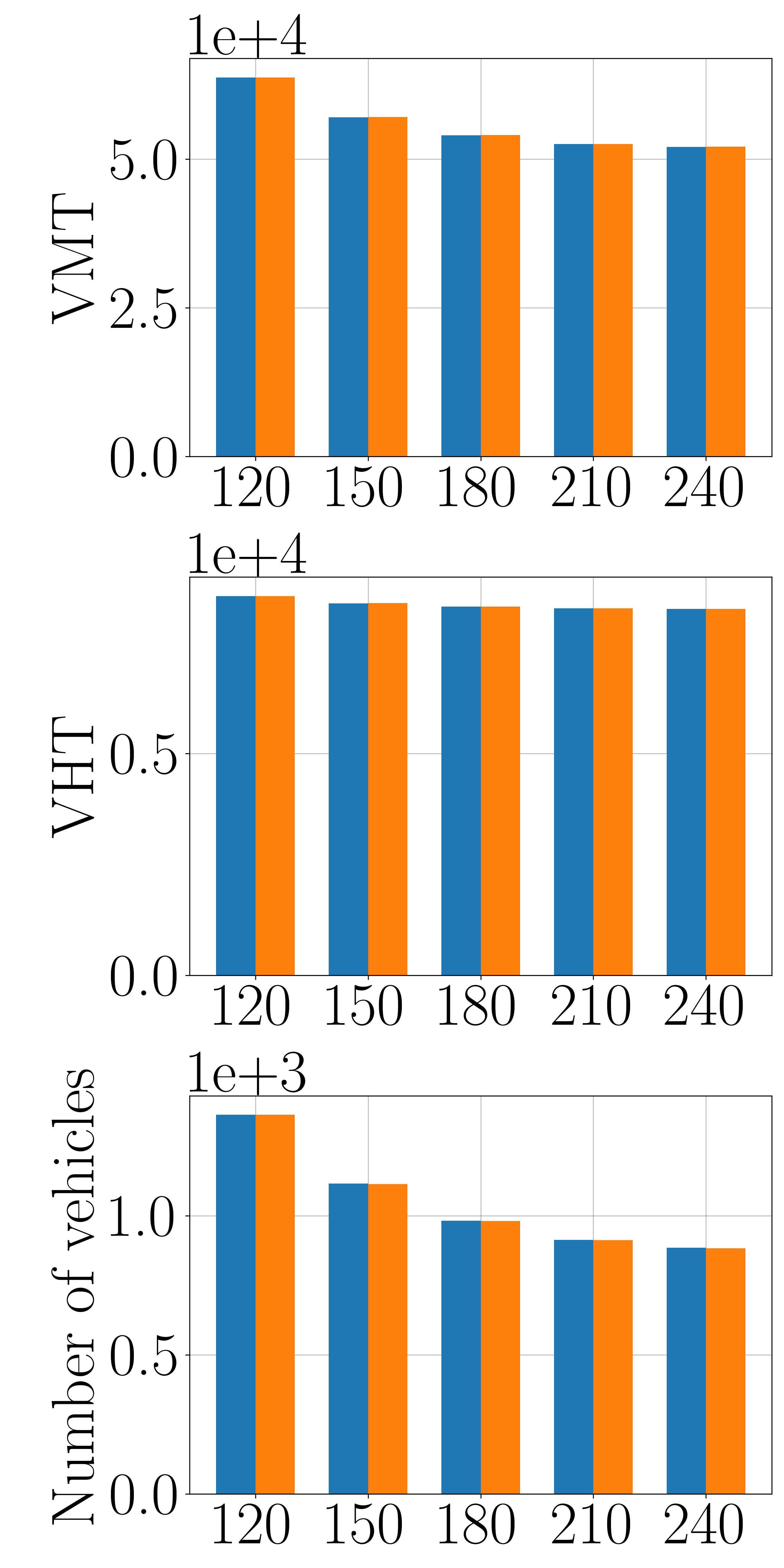}}
    \subfloat[\centering Bloomington\label{fig:Bloom_givenQ}]{%
        \includegraphics*[width=0.25\textwidth,height=\textheight,keepaspectratio]{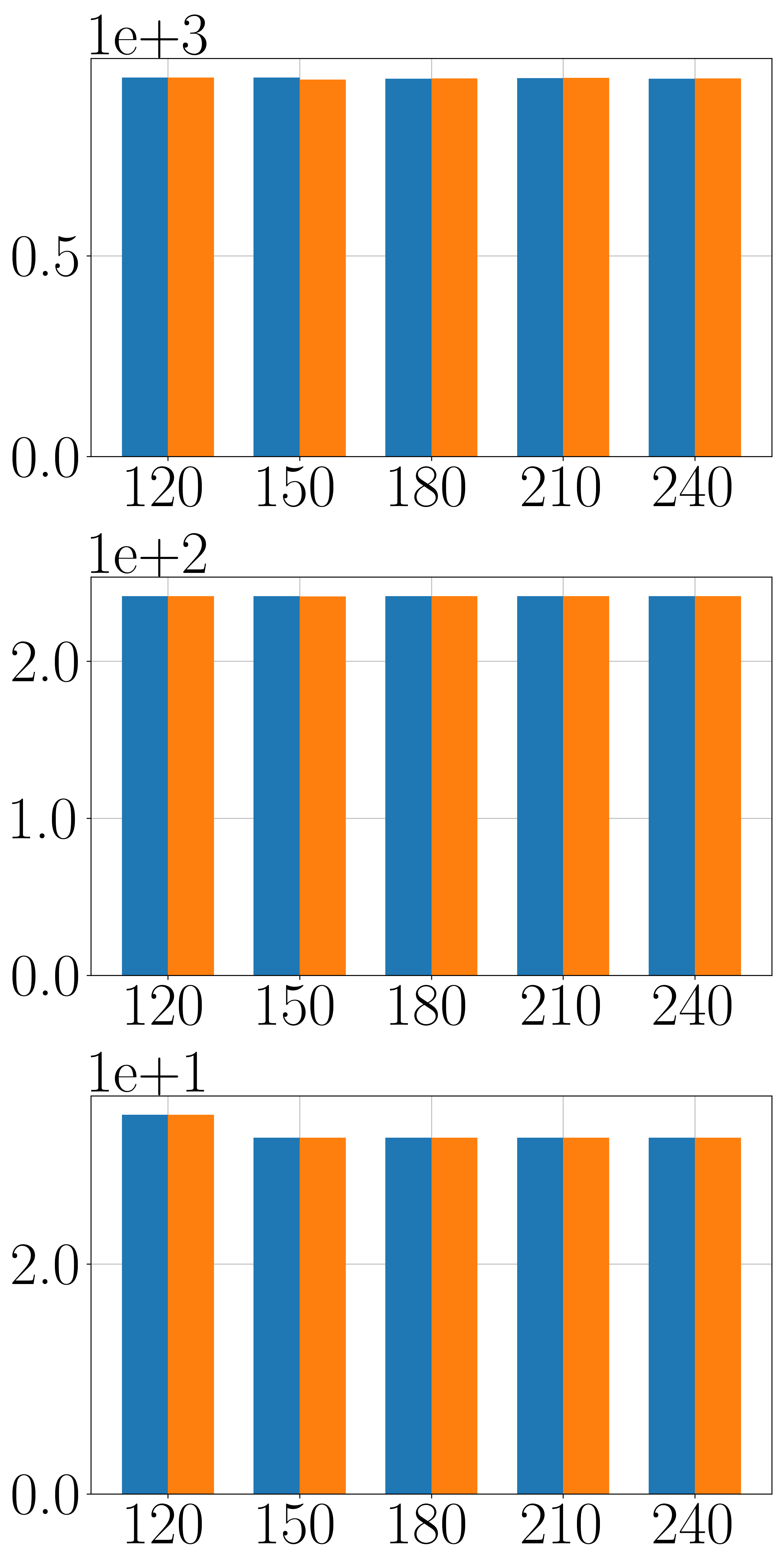}}
    \subfloat[\centering Chicago\label{fig:Chi_givenQ}]{%
        \includegraphics*[width=0.25\textwidth,height=\textheight,keepaspectratio]{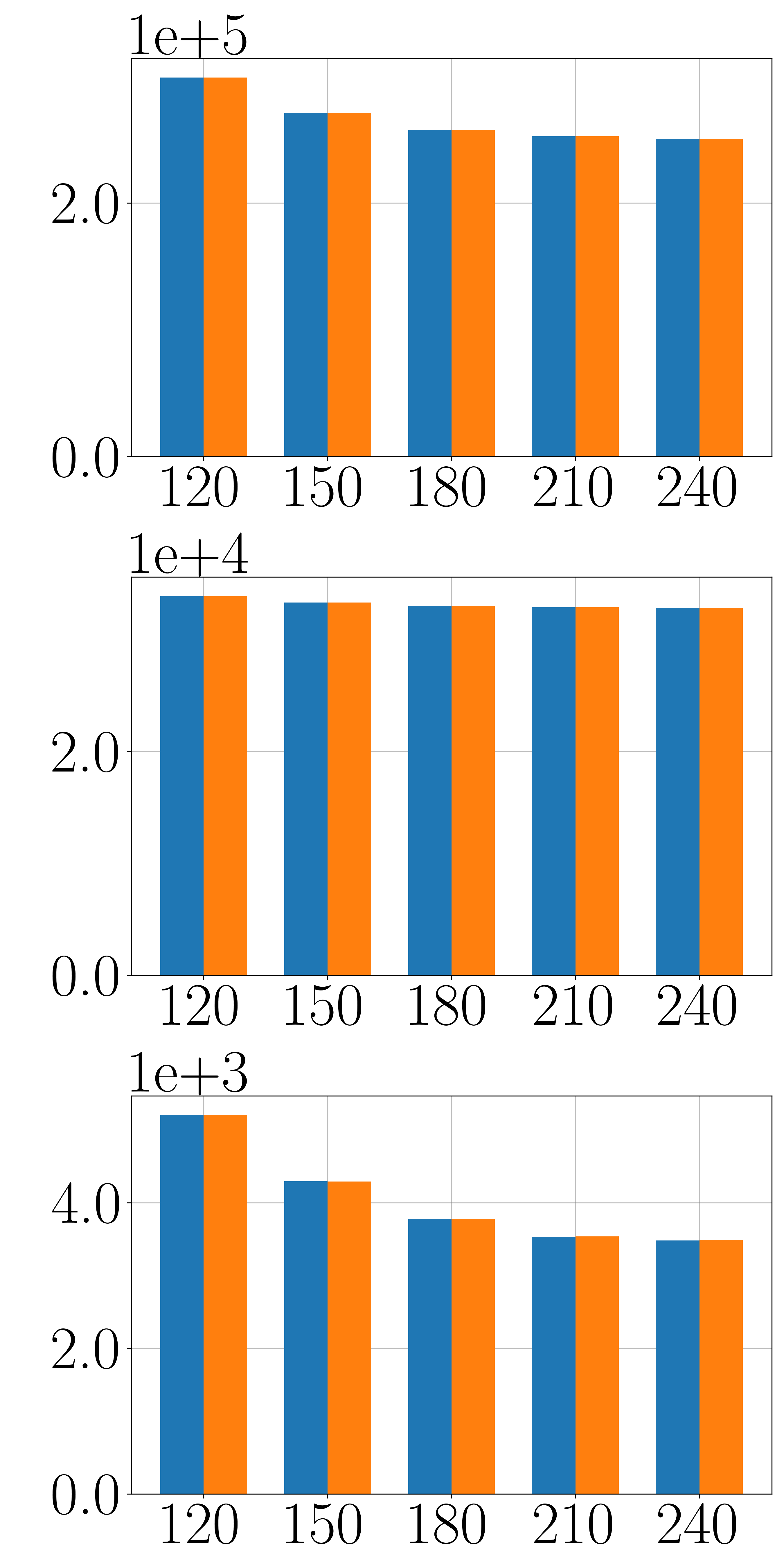}}
    \subfloat[\centering Detroit\label{fig:Det_givenQ}]{%
        \includegraphics*[width=0.25\textwidth,height=\textheight,keepaspectratio]{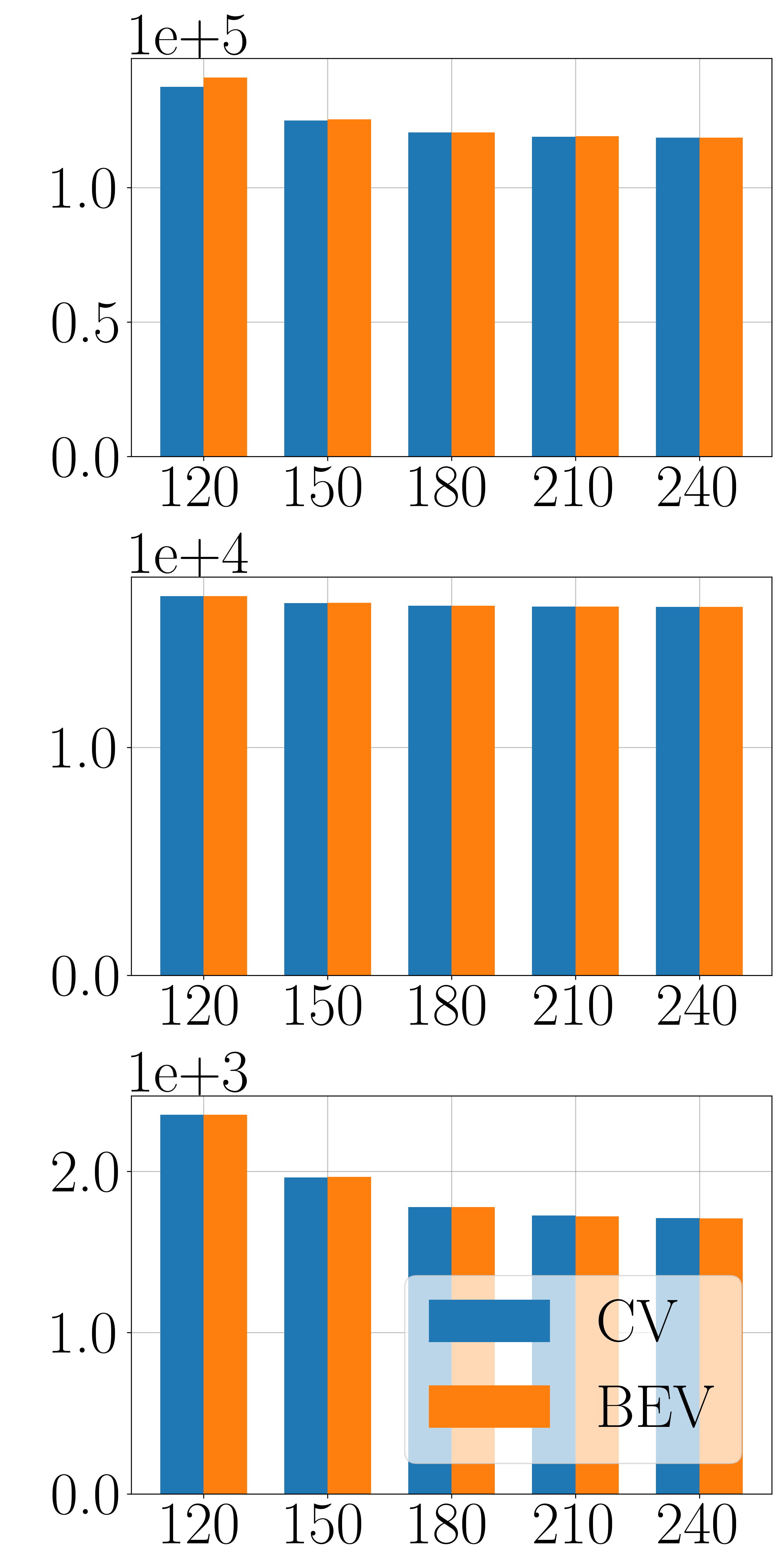}}
    \caption{Impact of vehicle capacity on VMT, VHT, and the number of customers.
    } \label[fig]{Results_givenQ}
\end{figure*}

Under the aforementioned experimental design, the vehicle capacity, $Q$, 
refers to the number of customers who can be served by a vehicle. 
We considered a set of vehicle capacity values, $Q\in\{120,~150,~180,~210,~240\}$, 
and solved all TCVRP instances in the four cities for BEVs and CVs. 
\cref{Results_givenQ} shows the system-level VMT, VHT, and the number of vehicles for these cases. 
The first impression is that BEV and CV metrics are almost the same in all cases. 
This is not surprising because the average VMT per vehicle,  VMT/(number of vehicles), is always below 80 miles. 
Therefore, BEV range constraints are not binding, and BEVs become equivalent to CVs. 
For example, in Chicago's BEV case with $Q=240$ (see \cref{fig:Chi_givenQ}), the average VMT per vehicle is 71.8 miles, which is also the maximum number across cities.

In \cref{fig:Bloom_givenQ} we observe that most metrics are unaffected by the capacity increase. 
The reason is that time constraints are binding for the majority of the vehicle routes when $Q=120$. 
Increasing $Q$ to 150 allows vehicles whose routes have binding capacity constraints 
(but nonbinding time constraints) at $Q=120$ to serve more customers. 
Therefore, the number of vehicles drops by 2, and it plateaus for $Q\geq 150$. Yet, 
such a decrease does not impact the overall VMT and VHT.

In the large cities, the results yield the expected impact of increased $Q$ 
(see \Cref{fig:Aus_givenQ,fig:Chi_givenQ,fig:Det_givenQ}), 
that is, (more or less) a decrease in all reported metrics. 
The system-level VMTs (and the energy consumption) in Austin, Chicago, and Detroit decrease by nearly 18, 16, and 14\%, 
respectively, when $Q$ doubles from 120 to 240. Similarly, 
fleet sizes decrease by 35, 33, and 27\% in the same order. Overall, we find that vehicle capacity constraints are binding for the majority of these instances.

\subsection{Impact of Maximum Allowed Travel Time on System-Level Metrics}
\begin{figure*}[!htbp]
    \centering
    \subfloat[\centering Austin\label{fig:Aus_givenT}]{%
        \includegraphics*[width=0.25\textwidth,height=\textheight,keepaspectratio]{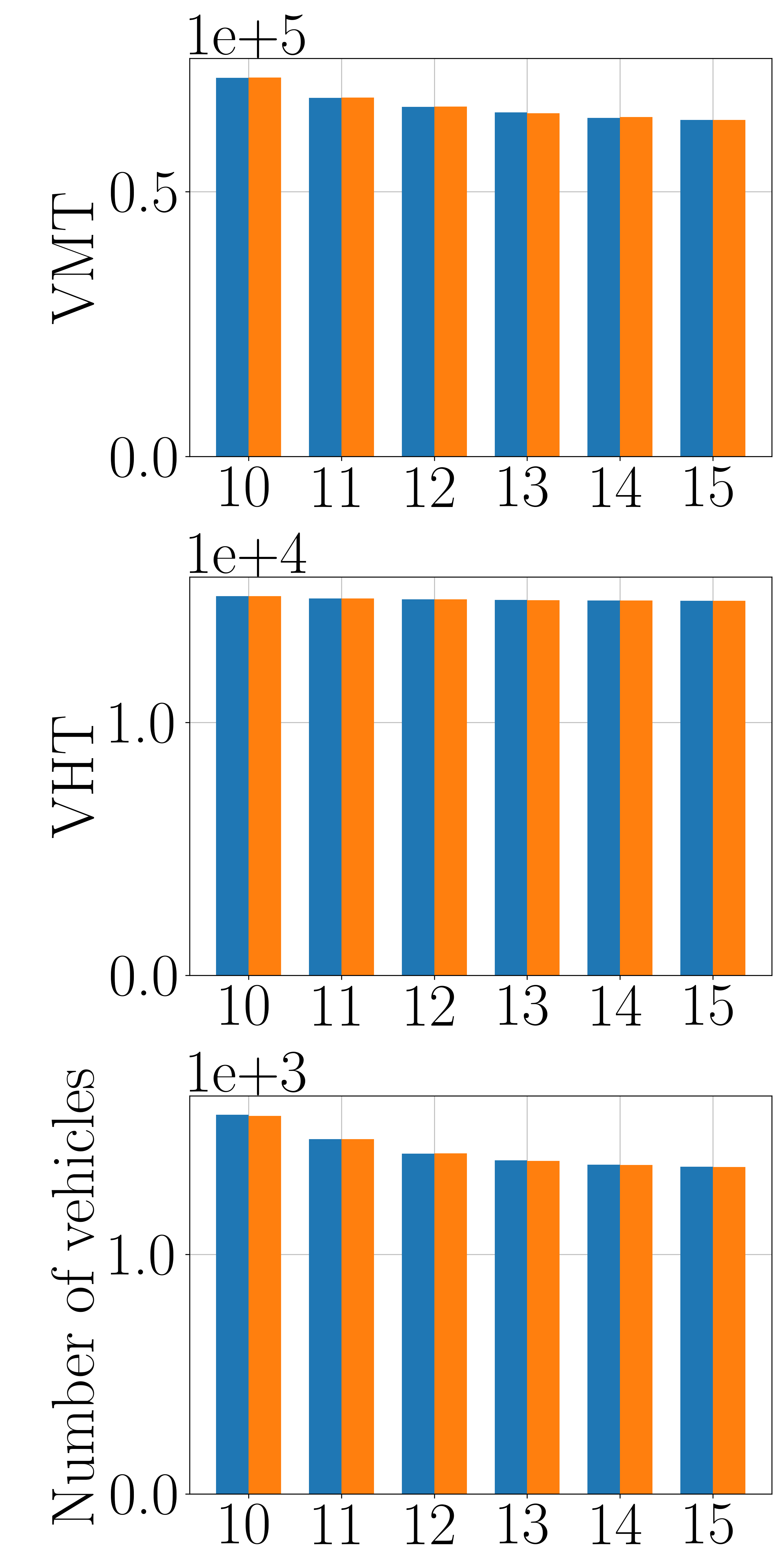}}
    \subfloat[\centering Bloomington\label{fig:Bloom_givenT}]{%
        \includegraphics*[width=0.25\textwidth,height=\textheight,keepaspectratio]{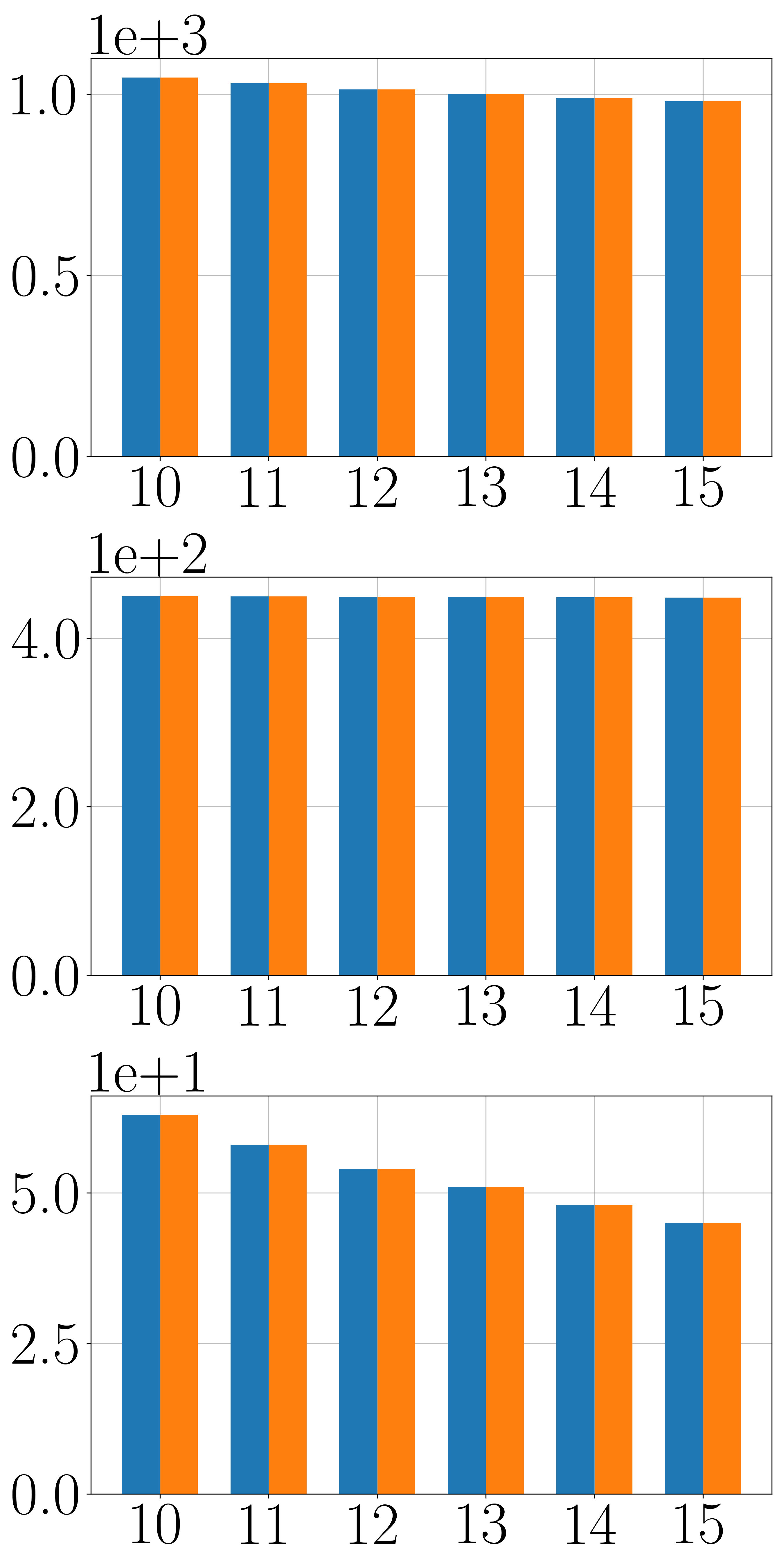}}
    \subfloat[\centering Chicago\label{fig:Chi_givenT}]{%
        \includegraphics*[width=0.25\textwidth,height=\textheight,keepaspectratio]{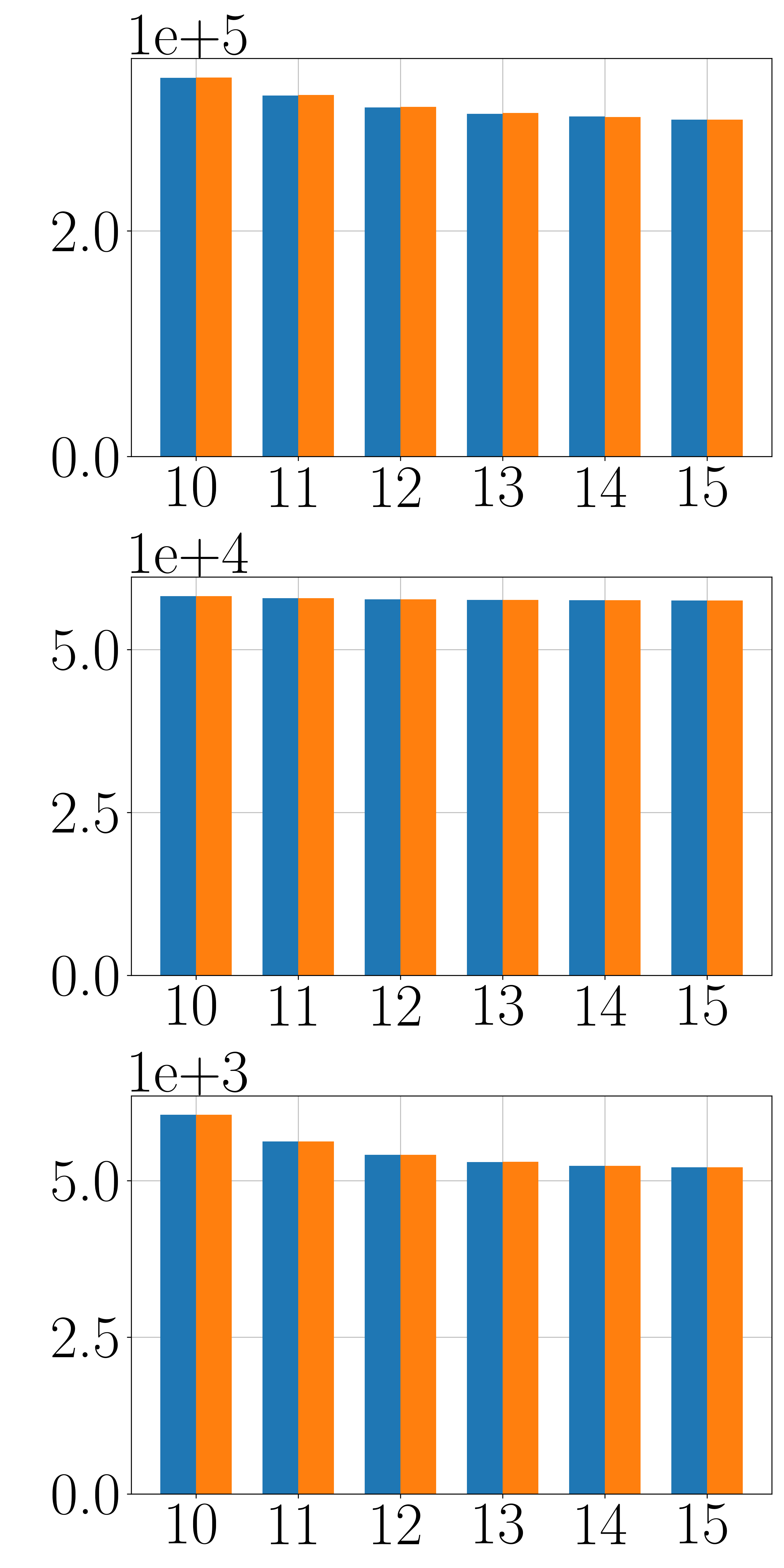}}
    \subfloat[\centering Detroit\label{fig:Det_givenT}]{%
        \includegraphics*[width=0.25\textwidth,height=\textheight,keepaspectratio]{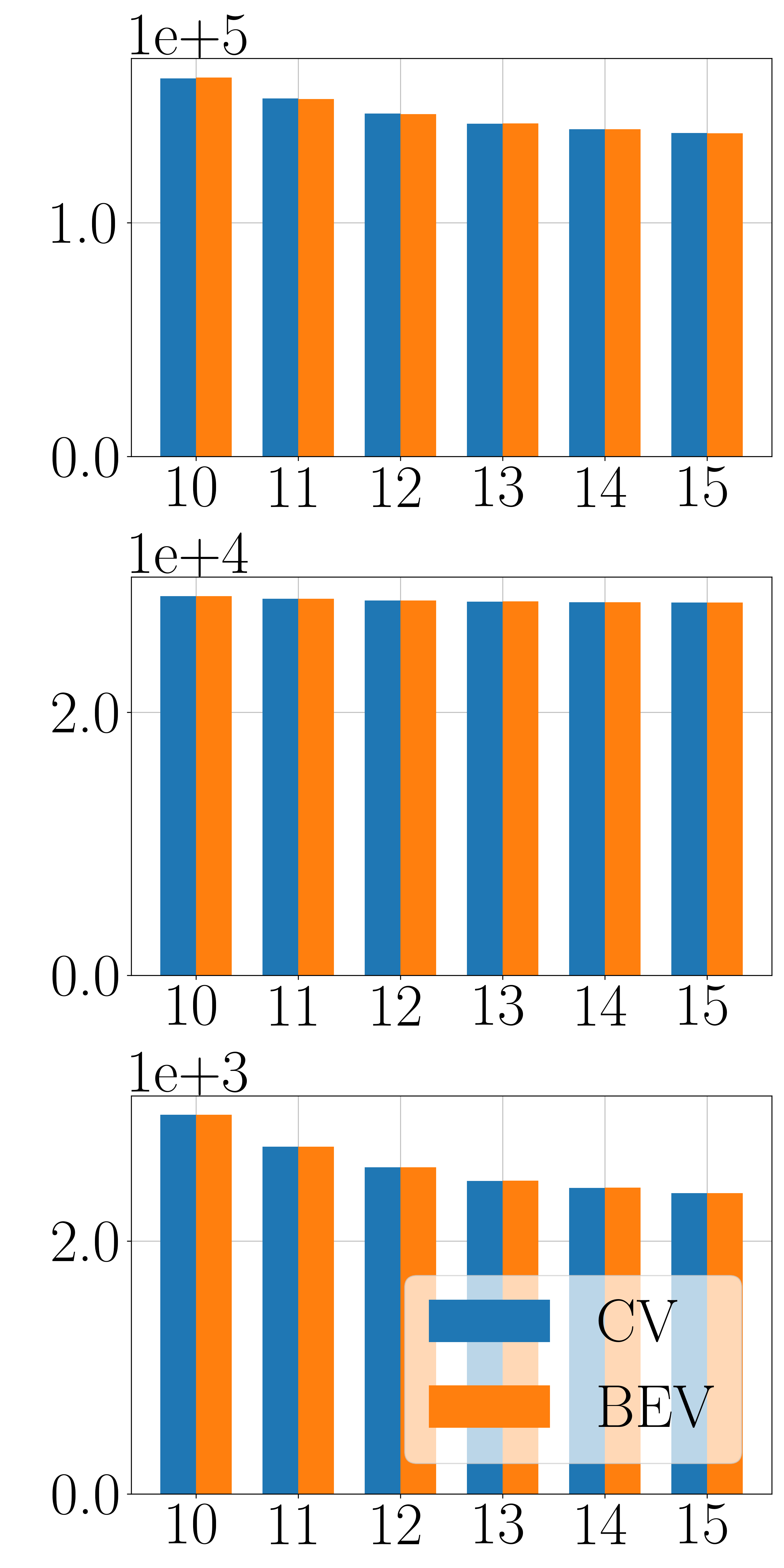}}
    \caption{Impact of maximum allowed travel time on VMT, VHT, and the number of customers.} \label[fig]{Results_givenT}
\end{figure*}

Maximum allowed travel time is a realistic constraint representing the limited work hours of service providers. 
Using the baseline parameter settings, we considered $\overline{T}\in\{10,~11,~12,~13,~14,~15\}$ hours in the four cities for both BEVs and CVs. 
The results, although not reported, show that  an increase does not impact the VMT, VHT, and the number of vehicles in the urban areas, 
and it has a minimal impact in the Bloomington case. This is because most of the vehicle routes have binding capacity constraints. 
To observe the expected impact, we changed $P$ from 2 minutes to 4 minutes in addition to testing the given values of $\overline{T}$. 
Once $P=4$ minutes and $\overline{T}=10$ hours, time constraints become the dominantly binding constraint; therefore, relaxing $\overline{T}$ yields the expected improvement in the reported key metrics. 
\cref{Results_givenT} shows the impact of increased $\overline{T}$ on these metrics. 
As expected, the VMT, VHT, and the number of vehicles decrease as $\overline{T}$ increases. 
VHT is the least impacted metric because it is an outcome of the minimized VMT, 
and a linear relation between the VMT and the VHT 
may not occur because arc travel speeds vary across the network.

\subsection{Impact of Service Time on System-Level Metrics}
\begin{figure*}[!htbp]
    \centering
    \subfloat[\centering Austin\label{fig:Aus_givenP}]{%
        \includegraphics*[width=0.25\textwidth,height=\textheight,keepaspectratio]{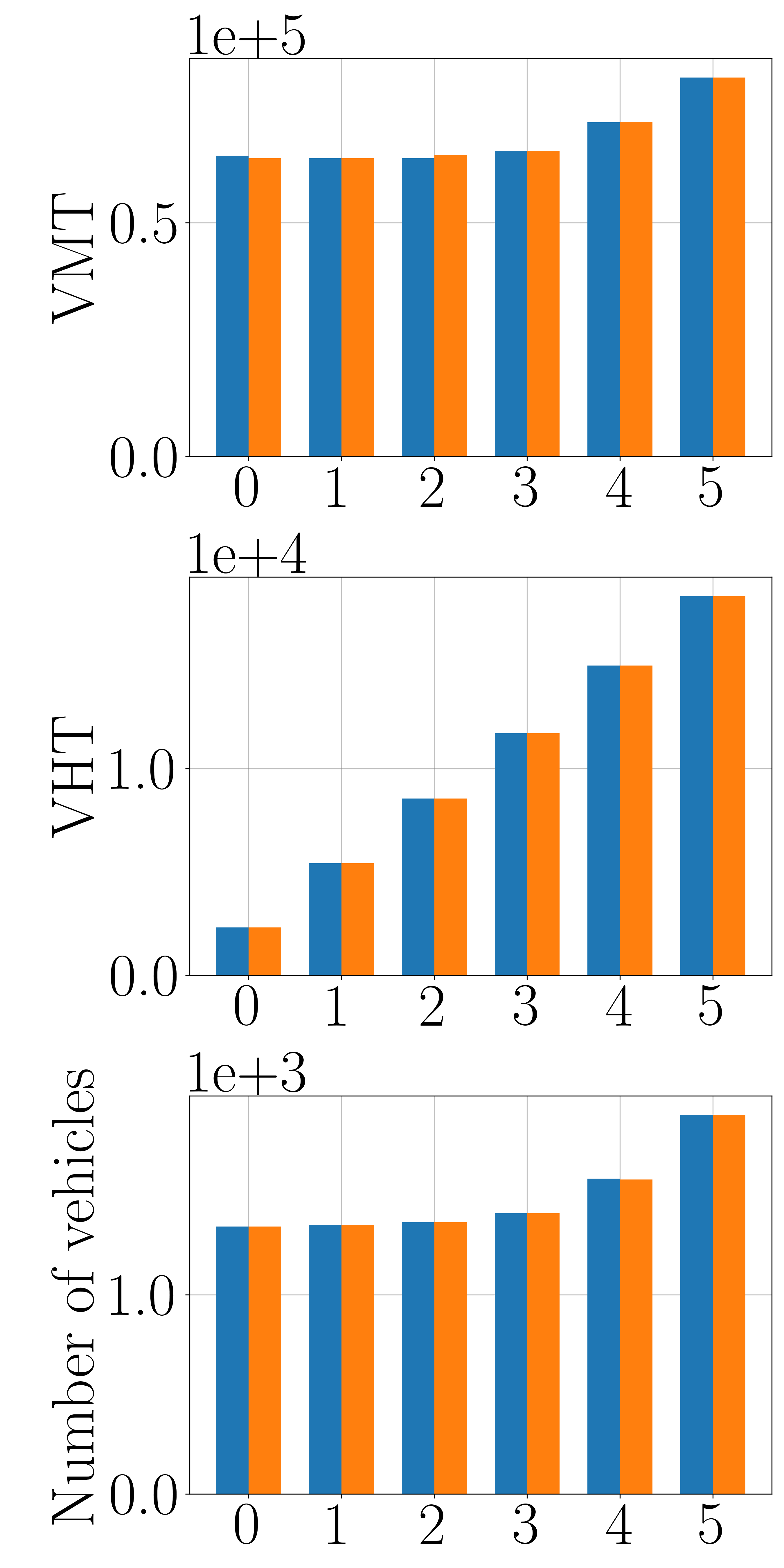}}
    \subfloat[\centering Bloomington\label{fig:Bloom_givenP}]{%
        \includegraphics*[width=0.25\textwidth,height=\textheight,keepaspectratio]{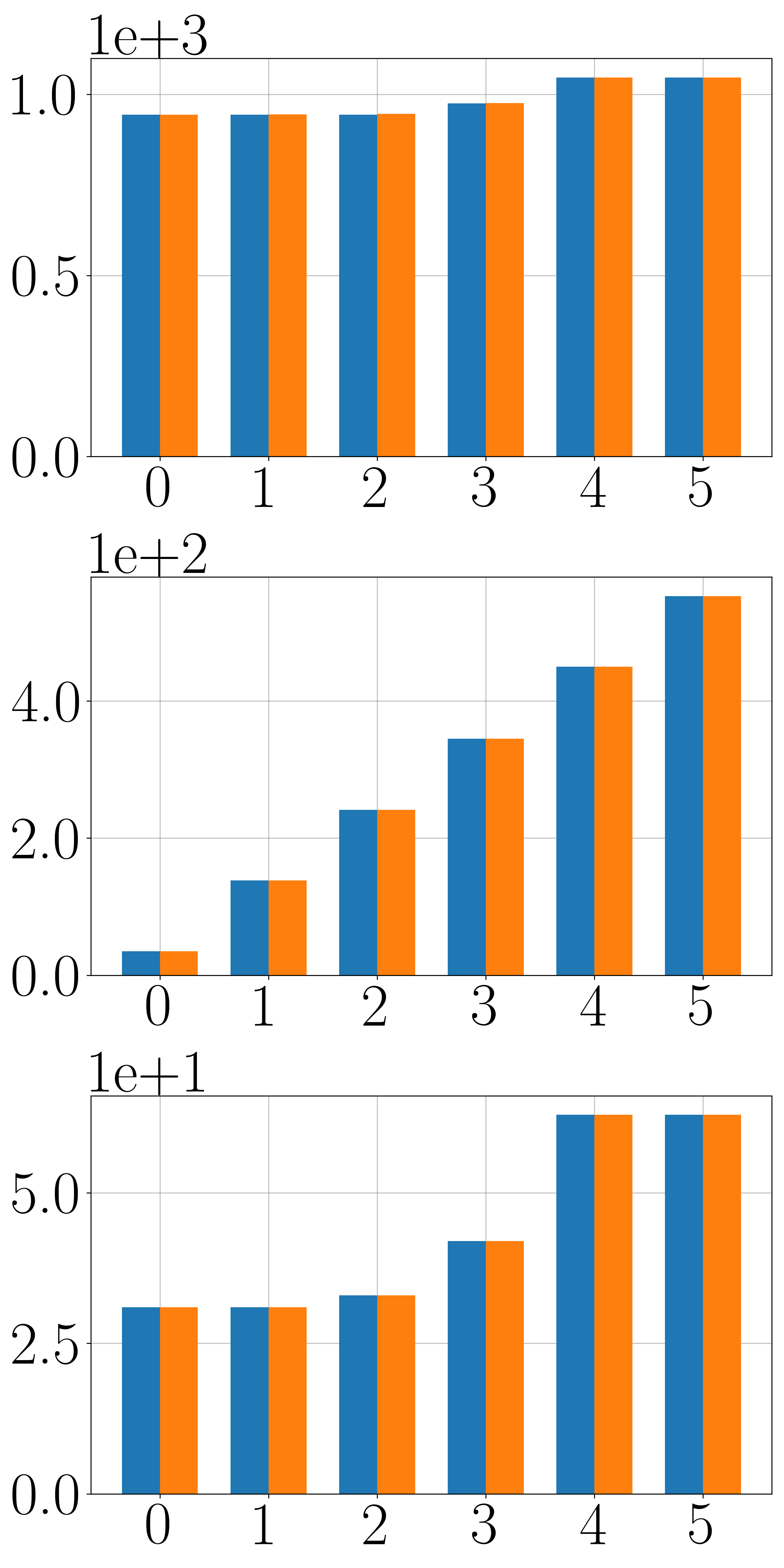}}
    \subfloat[\centering Chicago\label{fig:Chi_givenP}]{%
        \includegraphics*[width=0.25\textwidth,height=\textheight,keepaspectratio]{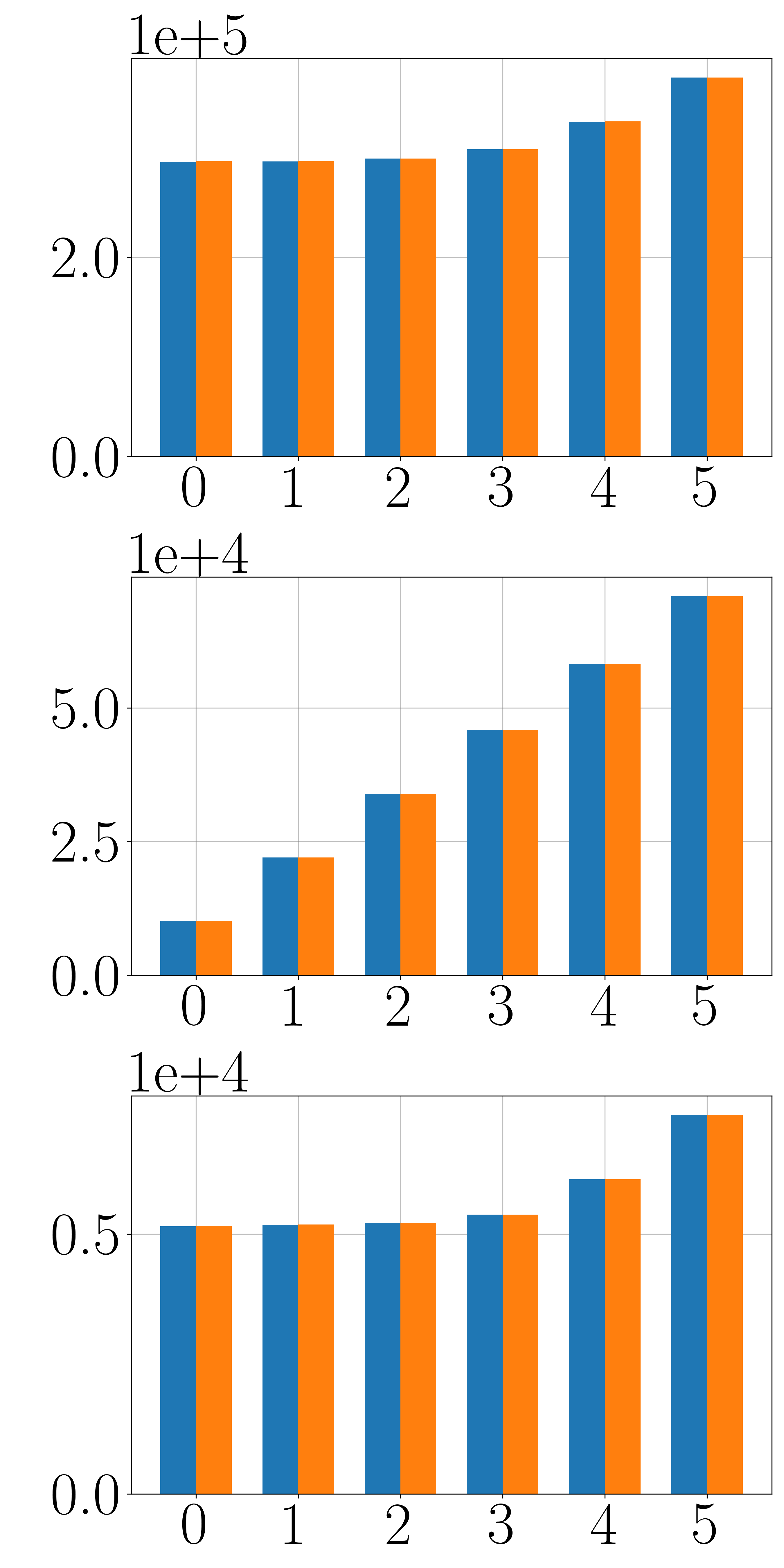}}
    \subfloat[\centering Detroit\label{fig:Det_givenP}]{%
        \includegraphics*[width=0.25\textwidth,height=\textheight,keepaspectratio]{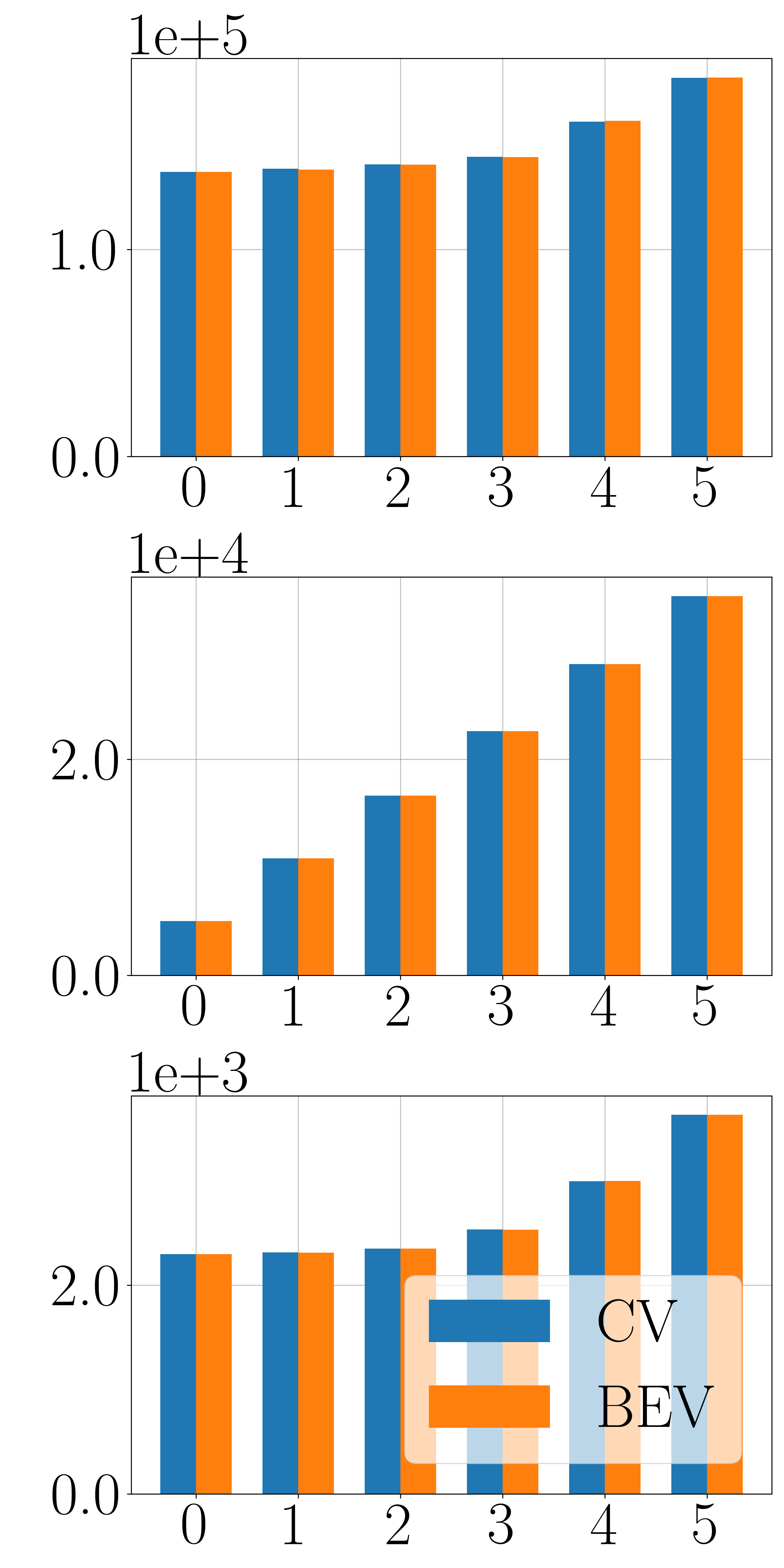}}
    \caption{Impact of service time on VMT, VHT, and the number of customers.} \label[fig]{Results_givenP}
\end{figure*}

Service time plays a critical role in the number of customers served by each vehicle. 
Assume that $P=5$ minutes and a vehicle's route includes 120 deliveries. 
Then, the corresponding service time is 10 hours, which forms an infeasible route in the current settings (i.e., $\overline{T}=10$ hours). 
For this reason, vehicles end up serving fewer customers when $P=5$. 
To investigate the impact of $P$ on system-level metrics, we considered $P\in\{0,~1,~2,~3,~4,~5\}$ minutes. 
\cref{Results_givenP} shows that $P$ has an exponentially increasing impact on the VMT and the number of vehicles, 
whereas it has a linear impact on the VHT 
because the service time 
per customer---increased linearly---is the dominant time factor in the travel time of vehicles; in  other words, 
the road travel time at the system level is far below the service time (especially when $P\geq 2$). 

The results in the figure show that five minutes of service time per customer (compared with $P=0$) increases 
the VMT by 27.1, 10.9, 28.4, and 33.1\% and the number of vehicles by 41.9, 103.2, 41.5, and 58.1\% in Austin, Bloomington, Chicago, and Detroit, respectively. 
These numbers indicate that urban areas are differently impacted by a change in $P$ compared with Bloomington. For example, a 103.2\% increase in the fleet size impacts the VMT by only 10.9\% which yields an approximate 10-to-1 ratio. The ratio in large cities, however, is around 1.5-to-1.

\subsection{Impact of BEV Range on System-Level Metrics}
Per the baseline parameter values and the earlier  sensitivity analyses, 
the highest systemwide average VMT/vehicle was around 72 miles in Chicago's case with $Q=240$. 
Most cases had binding  time or capacity constraints, and the impact of the BEV range constraints was limited. 
We considered longer BEV ranges than 80 miles and observed that the performance metrics remained almost unchanged.
In the large cities, we discovered the parameter values that make the VMT/vehicle near 80 miles. 
We did not consider the Bloomington case because the maximum systemwide VMT/vehicle was around 30 miles, 
and altering parameters within realistic boundaries would not yield a solution with a VMT/vehicle of 80 miles.
To maximize the VMT/vehicle, we should assume longer work hours, shorter service times, higher vehicle capacities,
or a mix of these assumptions. 
For the three large cities, we set $P=1$ minute per customer and kept $\overline{T}=10$ hours (as in the baseline). 
By testing different values of $Q$ for each city, we found that the VMT/vehicle approaches 80 miles 
when $Q=400$ in Austin, $Q=240$ in Chicago, and $Q=210$ in Detroit. 
These numbers reveal the fact that the BEV range becomes a dominant factor 
when vehicles possess such delivery capacities under the parametric design explored.

\subsection{Impact of Shared Economy on System-Level Metrics}
In a shared economy environment, providers can use depots of each other under predetermined conditions and pricing policies. 
Such shared use of resources could also bring benefits in the e-commerce parcel delivery context. 
Assume a provider can rely on another one to make deliveries for customers that are geographically closer to their depots. 
Then, the overall VMT would be expected to decrease. 
This system can also be considered as a centrally controlled parcel delivery system. 
We therefore analyzed the magnitude of reduction in VMT, VHT, and the number of vehicles when all deliveries were controlled by a central mechanism. 
We selected Austin as the case area, distributed customers to depots without differentiating the service providers, 
and  solved the resulting problem instance using the baseline parameter settings. 
Compared with its counterpart instance results, the centralization reduced the VMT and the VHT by 38.8\% and 19.7\%, respectively, and slightly (by 0.5\%) increased the number of vehicles.

\section{Conclusion}\label{sec:conclusion}
In this study, we developed an MIP as an exact solution method to solve the TCVRP of e-commerce parcel delivery BEVs and CVs.
We compared our method with a previously developed ITS metaheuristic and presented the performance statistics. 
Although the ITS performed better in most instances, 
the MIP was found useful to prove optimality and ensure a confidence level for the solutions obtained by the ITS. 
Supported by validated simulation data of POLARIS, we designed an experimental layout and analyzed three large cities---Austin, Chicago, and Detroit---and the smaller city of Bloomington. 
Because of large problem sizes, we aggregated customers at arc midpoints called 
super-locations and solved all problem instances using the ITS metaheuristic. 

We considered the impact of vehicle capacity, maximum allowed travel time, service time, 
and BEV range on the system-level metrics (i.e., VMT, VHT, and the number of vehicles). 
The results in the four cities showed that the service time followed by the vehicle 
capacity impacts the system-level metrics the most. In Bloomington, increasing vehicle capacities did not impact the system-level metrics while other areas have shown considerable reduction in VMT, VHT, and the number of vehicles. On the other hand, the increase in the dwelling time impacted Bloomington the most compared to the other cities. This can be related to low customer density in the area and the size of the city.

Simplifying the EVRP by omitting recharging decisions, 
we solved the BEV routing problem under service time, capacity, and BEV range constraints. 
Case studies illustrated that the BEV range is not a limiting factor since the maximum of 
the average VMT per vehicle across scenarios was around 72 miles. 
Yet, we also identified that increased vehicle capacities with a 
dwelling time of one minute (per customer) can alter the picture and make 
the BEV range an important factor that impacts system-level metrics. 
Based on our findings, vehicle capacities of 400, 240, and 210 for Austin, Chicago, and Detroit, respectively, 
are the breakpoints where the BEV range begins gaining importance. 
Technically, most BEV range constraints become binding when each vehicle has a higher capacity than the one denoted for each area.

We extended the analyses by accounting for a shared economy scenario considering Austin as the case study area. 
In this scenario, we assumed that all deliveries can be made by any depot in the region regardless of the ownership of the depot. Although this analysis is not relevant to model parameters, it yields a significant reduction in the key metrics that is worth sharing.
Under these assumptions, the results indicated that VMT and VHT decrease by 38.8\% and 19.7\%, respectively.

We provided an energy consumption estimate based on the multiplication of the VMT and the kWh per mile energy consumption of BEVs. 
Also, we denoted the CV fuel efficiency and the energy unit equivalence of the diesel to account for an approximate energy consumption of CVs. 
For each city, we identified the average VMT across the scenarios considered and found the following: 
BEVs consume 0.07, 0.001, 0.34, and 0.16 gWh, and CVs consume 0.32, 0.005, 1.51, and 0.71 gWh in Austin, Bloomington, Chicago, and Detroit, respectively.

Multiple possibilities remain for enhancement in the modeling and analyses. 
The exact method can certainly be improved by considering a route-based formulation as in~\cite{subramanyam2021joint}. 
A more comprehensive analysis could account for varying service times and other parameters simultaneously 
to better observe the impact of these parameters on the system-level performance metrics. 
The results showed that different outcomes can be observed based on the areas considered. 
For instance, BEVs may not be useful in small cities such as Bloomington because the average VMT per vehicle is low.
Hence, extending these analyses to other cities would better inform decision-makers about areas where BEVs could be suitable.
Our findings can be interpreted as follows: Urban areas can benefit from BEVs more since the VMT per vehicle is higher compared with that of small cities. Yet, we should also state the fact that VMT per vehicle could also be higher in small cities depending on the spatial distribution of customer and depots as well as their density.

\section{Acknowledgements}
This material is based upon work supported by the U.S. Department of Energy, Office of Science, under contract number DE-AC02-06CH11357. 
This report and the work described were sponsored by the U.S. Department of Energy (DOE) Vehicle Technologies Office (VTO) under the Systems and Modeling for Accelerated Research in Transportation (SMART) Mobility Laboratory Consortium, an initiative of the Energy Efficient Mobility Systems (EEMS) Program. 
Erin Boyd, a DOE Office of Energy Efficiency and Renewable Energy (EERE) manager, played an important role in establishing the project concept, advancing implementation, and providing guidance.

\section{Author Contributions}
The authors confirm contribution to the paper as follows: Study conception and
design: T. Cokyasar and M. Stinson; data collection: T. Cokyasar, O. Sahin, and M.
Stinson; analysis and interpretation of results: T. Cokyasar and A. Subramanyam; draft manuscript
preparation: T. Cokyasar and J. Larson. All authors reviewed the results and
approved the final version of the manuscript.

\section*{ORCID}

\def\orcid#1{\kern .08em\href{https://orcid.org/#1}{\includegraphics[keepaspectratio,width=0.7em]{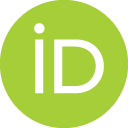}}}
Taner Cokyasar \orcid{0000-0001-9687-6725}, Anirudh Subramanyam \orcid{0000-0001-5255-8821}, Jeffrey Larson \orcid{0000-0001-9924-2082}, Monique Stinson \orcid{0000-0003-1337-1903}, Olcay Sahin \orcid{0000-0003-4235-036X}

\bibliographystyle{TRR}
\bibliography{TRR_references}

\begin{thebibliography}{10}
\providecommand{\url}[1]{\texttt{#1}}
\providecommand{\urlprefix}{URL }
\expandafter\ifx\csname urlstyle\endcsname\relax
  \providecommand{\doi}[1]{doi:\discretionary{}{}{}#1}\else
  \providecommand{\doi}{doi:\discretionary{}{}{}\begingroup
  \urlstyle{rm}\Url}\fi
\providecommand{\eprint}[2][]{\url{#2}}

\bibitem{archetti2011complexity}
Archetti, C., D.~Feillet, M.~Gendreau, and M.~G. Speranza.
\newblock Complexity of the {VRP} and {SDVRP}.
\newblock \emph{Transportation Research Part C: Emerging Technologies},
  Vol.~19, No.~5, 2011, pp. 741--750.
\newblock \doi{10.1016/j.trc.2009.12.006}.

\bibitem{almoustafa2013new}
Almoustafa, S., S.~Hanafi, and N.~Mladenovi{\'c}.
\newblock New exact method for large asymmetric distance-constrained vehicle
  routing problem.
\newblock \emph{European Journal of Operational Research}, Vol. 226, No.~3,
  2013, pp. 386--394.
\newblock \doi{10.1016/j.ejor.2012.11.040}.

\bibitem{kara2007energy}
Kara, I., B.~Y. Kara, and M.~K. Yetis.
\newblock Energy minimizing vehicle routing problem.
\newblock In \emph{International Conference on Combinatorial Optimization and
  Applications}. Springer, 2007, pp. 62--71.
\newblock \doi{10.1007/978-3-540-73556-4_9}.

\bibitem{kara2011arc}
Kara, I.
\newblock Arc based integer programming formulations for the distance
  constrained vehicle routing problem.
\newblock In \emph{International Symposium on Logistics and Industrial
  informatics}. IEEE, 2011, pp. 33--38.
\newblock \doi{10.1109/lindi.2011.6031159}.

\bibitem{kek2008distance}
Kek, A.~G., R.~L. Cheu, and Q.~Meng.
\newblock Distance-constrained capacitated vehicle routing problems with
  flexible assignment of start and end depots.
\newblock \emph{Mathematical and Computer Modelling}, Vol.~47, No. 1-2, 2008,
  pp. 140--152.
\newblock \doi{10.1016/j.mcm.2007.02.007}.

\bibitem{arnold2019efficiently}
Arnold, F., M.~Gendreau, and K.~S{\"o}rensen.
\newblock Efficiently solving very large-scale routing problems.
\newblock \emph{Computers \& Operations Research}, Vol. 107, 2019, pp. 32--42.
\newblock \doi{10.1016/j.cor.2019.03.006}.

\bibitem{auld2016polaris}
Auld, J., M.~Hope, H.~Ley, V.~Sokolov, B.~Xu, and K.~Zhang.
\newblock {POLARIS: A}gent-based Modeling Framework Development and
  Implementation for Integrated Travel Demand and Network and Operations
  Simulations.
\newblock \emph{Transportation Research Part C: Emerging Technologies},
  Vol.~64, 2016, pp. 101--116.
\newblock \doi{10.1016/j.trc.2015.07.017}.

\bibitem{dantzig1959truck}
Dantzig, G.~B. and J.~H. Ramser.
\newblock The truck dispatching problem.
\newblock \emph{Management Science}, Vol.~6, No.~1, 1959, pp. 80--91.
\newblock \doi{10.1287/mnsc.6.1.80}.

\bibitem{flood1956traveling}
Flood, M.~M.
\newblock The traveling-salesman problem.
\newblock \emph{Operations Research}, Vol.~4, No.~1, 1956, pp. 61--75.
\newblock \doi{10.1287/opre.4.1.61}.

\bibitem{bertsimas1992vehicle}
Bertsimas, D.~J.
\newblock A vehicle routing problem with stochastic demand.
\newblock \emph{Operations Research}, Vol.~40, No.~3, 1992, pp. 574--585.
\newblock \doi{10.1287/opre.40.3.574}.

\bibitem{ordonez2007priori}
Ord{\'{o}}{\~{n}}ez, F., I.~Sungur, and M.~Dessouky.
\newblock A priori performance measures for arc-based formulations of vehicle
  routing problem.
\newblock \emph{Transportation Research Record}, Vol. 2032, No.~1, 2007, pp.
  53--62.
\newblock \doi{10.3141/2032-07}.

\bibitem{ehmke2018optimizing}
Ehmke, J.~F., A.~M. Campbell, and B.~W. Thomas.
\newblock Optimizing for total costs in vehicle routing in urban areas.
\newblock \emph{Transportation Research Part E: Logistics and Transportation
  Review}, Vol. 116, 2018, pp. 242--265.
\newblock \doi{10.1016/j.tre.2018.06.008}.

\bibitem{figliozzi2010vehicle}
Figliozzi, M.
\newblock Vehicle routing problem for emissions minimization.
\newblock \emph{Transportation Research Record}, Vol. 2197, No.~1, 2010, pp.
  1--7.
\newblock \doi{10.3141/2197-01}.

\bibitem{drexl2012synchronization}
Drexl, M.
\newblock Synchronization in Vehicle Routing{\textemdash}A Survey of {VRPs}
  with Multiple Synchronization Constraints.
\newblock \emph{Transportation Science}, Vol.~46, No.~3, 2012, pp. 297--316.
\newblock \doi{10.1287/trsc.1110.0400}.

\bibitem{wang2013vehicle}
Wang, Y., X.~Ma, Y.~Lao, Y.~Wang, and H.~Mao.
\newblock Vehicle routing problem: Simultaneous deliveries and pickups with
  split loads and time windows.
\newblock \emph{Transportation Research Record}, Vol. 2378, No.~1, 2013, pp.
  120--128.
\newblock \doi{10.3141/2378-13}.

\bibitem{eksioglu2009vehicle}
Eksioglu, B., A.~V. Vural, and A.~Reisman.
\newblock The vehicle routing problem: A taxonomic review.
\newblock \emph{Computers \& Industrial Engineering}, Vol.~57, No.~4, 2009, pp.
  1472--1483.
\newblock \doi{10.1016/j.cie.2009.05.009}.

\bibitem{braekers2016vehicle}
Braekers, K., K.~Ramaekers, and I.~Van~Nieuwenhuyse.
\newblock The vehicle routing problem: State of the art classification and
  review.
\newblock \emph{Computers \& Industrial Engineering}, Vol.~99, 2016, pp.
  300--313.
\newblock \doi{10.1016/j.cie.2015.12.007}.

\bibitem{konstantakopoulos2020vehicle}
Konstantakopoulos, G.~D., S.~P. Gayialis, and E.~P. Kechagias.
\newblock Vehicle routing problem and related algorithms for logistics
  distribution: A literature review and classification.
\newblock \emph{Operational Research}, 2020, pp. 1--30.
\newblock \doi{10.1007/s12351-020-00600-7}.

\bibitem{laporte1987branch}
Laporte, G., Y.~Nobert, and S.~Taillefer.
\newblock A branch-and-bound algorithm for the asymmetrical
  distance-constrained vehicle routing problem.
\newblock \emph{Mathematical Modelling}, Vol.~9, No.~12, 1987, pp. 857--868.
\newblock \doi{10.1016/0270-0255(87)90004-2}.

\bibitem{hashimoto2006vehicle}
Hashimoto, H., T.~Ibaraki, S.~Imahori, and M.~Yagiura.
\newblock The vehicle routing problem with flexible time windows and traveling
  times.
\newblock \emph{Discrete Applied Mathematics}, Vol. 154, No.~16, 2006, pp.
  2271--2290.
\newblock \doi{10.1016/j.dam.2006.04.009}.

\bibitem{kara2011polynomial}
Kara, I. and T.~Derya.
\newblock Polynomial size formulations for the distance and capacity
  constrained vehicle routing problem.
\newblock In \emph{AIP Conference Proceedings}, Vol. 1389. American Institute
  of Physics, 2011, pp. 1713--1718.
\newblock \doi{10.1063/1.3636940}.

\bibitem{braysy2005vehicle}
Br{\"a}ysy, O. and M.~Gendreau.
\newblock Vehicle routing problem with time windows, Part I: Route construction
  and local search algorithms.
\newblock \emph{Transportation science}, Vol.~39, No.~1, 2005, pp. 104--118.
\newblock \doi{10.1287/trsc.1030.0056}.

\bibitem{braysy2005vehicle2}
Br{\"a}ysy, O. and M.~Gendreau.
\newblock Vehicle routing problem with time windows, Part II: Metaheuristics.
\newblock \emph{Transportation science}, Vol.~39, No.~1, 2005, pp. 119--139.
\newblock \doi{10.1287/trsc.1030.0057}.

\bibitem{kallehauge2008formulations}
Kallehauge, B.
\newblock Formulations and exact algorithms for the vehicle routing problem
  with time windows.
\newblock \emph{Computers \& Operations Research}, Vol.~35, No.~7, 2008, pp.
  2307--2330.
\newblock \doi{10.1016/j.cor.2006.11.006}.

\bibitem{chen2016electric}
Chen, J., M.~Qi, and L.~Miao.
\newblock The electric vehicle routing problem with time windows and battery
  swapping stations.
\newblock In \emph{International Conference on Industrial Engineering and
  Engineering Management}. IEEE, 2016, pp. 712--716.
\newblock \doi{10.1109/IEEM.2016.7797968}.

\bibitem{keskin2018matheuristic}
Keskin, M. and B.~{\c{C}}atay.
\newblock A matheuristic method for the electric vehicle routing problem with
  time windows and fast chargers.
\newblock \emph{Computers \& Operations Research}, Vol. 100, 2018, pp.
  172--188.
\newblock \doi{10.1016/j.cor.2018.06.019}.

\bibitem{loffler2020routing}
L{\"o}ffler, M., G.~Desaulniers, S.~Irnich, and M.~Schneider.
\newblock Routing electric vehicles with a single recharge per route.
\newblock \emph{Networks}, Vol.~76, No.~2, 2020, pp. 187--205.
\newblock \doi{10.1002/net.21964}.

\bibitem{barnitt2011fedex}
Barnitt, R.
\newblock \emph{Fedex express gasoline hybrid electric delivery truck
  evaluation: 12-month report}.
\newblock Tech. rep., National Renewable Energy Lab, 2011.
\newblock \doi{10.2172/1007342}.

\bibitem{feng2013economic}
Feng, W. and M.~Figliozzi.
\newblock An economic and technological analysis of the key factors affecting
  the competitiveness of electric commercial vehicles: A case study from the
  USA market.
\newblock \emph{Transportation Research Part C: Emerging Technologies},
  Vol.~26, 2013, pp. 135--145.
\newblock \doi{10.1016/j.trc.2012.06.007}.

\bibitem{davis2013methodology}
Davis, B.~A. and M.~A. Figliozzi.
\newblock A methodology to evaluate the competitiveness of electric delivery
  trucks.
\newblock \emph{Transportation Research Part E: Logistics and Transportation
  Review}, Vol.~49, No.~1, 2013, pp. 8--23.
\newblock \doi{10.1016/j.tre.2012.07.003}.

\bibitem{kucukoglu2021electric}
Kucukoglu, I., R.~Dewil, and D.~Cattrysse.
\newblock The electric vehicle routing problem and its variations: A literature
  review.
\newblock \emph{Computers \& Industrial Engineering}, 2021, p. 107650.
\newblock \doi{10.1016/j.cie.2021.107650}.

\bibitem{chen2014preliminary}
Chen, Q. and J.~Lin.
\newblock A preliminary investigation of sustainable urban truck routing
  strategies considering cargo weight and vehicle speed.
\newblock In \emph{Transportation Research Board Annual Meeting}. 14-3300,
  2014, pp. 1--18.

\bibitem{nhts}
{Federal Highway Administration}.
\newblock National Household Travel Survey, 2017.
\newblock Available at \url{https://nhts.ornl.gov}, accessed on Oct. 28, 2021.

\bibitem{pitney}
Spadafora, J. and M.~Rodriguez.
\newblock Pitney Bowes Parcel Shipping Index Reveals 37 Percent Parcel Volume
  Growth in US for 2020, 2021.
\newblock Available at
  \url{https://www.businesswire.com/news/home/20210914005274/en/Pitney-Bowes-Parcel-Shipping-Index-Reveals-37-Percent-Parcel-Volume-Growth-in-US-for-2020},
  accessed on Dec. 10, 2021.

\bibitem{giosa2002new}
Giosa, I., I.~Tansini, and I.~Viera.
\newblock New assignment algorithms for the multi-depot vehicle routing
  problem.
\newblock \emph{Journal of the operational research society}, Vol.~53, No.~9,
  2002, pp. 977--984.
\newblock \doi{10.1057/palgrave.jors.2601426}.

\bibitem{drexl2015survey}
Drexl, M. and M.~Schneider.
\newblock A survey of variants and extensions of the location-routing problem.
\newblock \emph{European Journal of Operational Research}, Vol. 241, No.~2,
  2015, pp. 283--308.
\newblock \doi{10.1016/j.ejor.2014.08.030}.

\bibitem{laporte1986exact}
Laporte, G., Y.~Nobert, and D.~Arpin.
\newblock An exact algorithm for solving a capacitated location-routing
  problem.
\newblock \emph{Annals of Operations Research}, Vol.~6, No.~9, 1986, pp.
  291--310.
\newblock \doi{10.1007/BF02023807}.

\bibitem{pythontsp}
Goulart, F.
\newblock {Simple library to solve the Traveling Salesperson Problem in pure
  Python}, 2021.
\newblock Available at \url{https://github.com/fillipe-gsm/python-tsp} accessed
  on Dec. 24, 2021.

\bibitem{gurobi_tsp}
Gurobi.
\newblock tsp.py, 2021.
\newblock Available at
  \url{https://www.gurobi.com/documentation/9.1/examples/tsp_py.html} accessed
  on Jun. 2, 2021.

\bibitem{Subramanyam2020}
Subramanyam, A., P.~P. Repoussis, and C.~E. Gounaris.
\newblock Robust Optimization of a Broad Class of Heterogeneous Vehicle Routing
  Problems Under Demand Uncertainty.
\newblock \emph{{INFORMS} Journal on Computing}, Vol.~32, No.~3, 2020, pp.
  661--681.
\newblock \doi{10.1287/ijoc.2019.0923}.

\bibitem{krok}
Krok, A.
\newblock UPS to deploy 50 plug-in hybrid delivery trucks, 2018.
\newblock Available at
  \url{https://www.cnet.com/roadshow/news/ups-plug-in-hybrid-delivery-trucks-workhorse-group/},
  accessed on Dec. 24, 2021.

\bibitem{eia}
{Energy Information Administration}.
\newblock Annual Energy Outlook 2021.
\newblock Available at \url{https://www.eia.gov/outlooks/aeo/} accessed on Oct.
  31, 2021, 2021.

\bibitem{eia2}
{Energy Information Administration}.
\newblock Units and calculators explained.
\newblock Available at
  \url{https://www.eia.gov/energyexplained/units-and-calculators/} accessed on
  Oct. 31, 2021, 2021.

\bibitem{cerc}
{U.S.-China Clean Energy Research Center}.
\newblock {Truck Research Utilizing Collaborative Knowledge}, 2021.
\newblock Available at \url{https://cerc-truck.anl.gov/} accessed on Dec. 23,
  2021.

\bibitem{gurobi}
{Gurobi Optimization, LLC}.
\newblock Gurobi Optimizer Reference Manual.
\newblock Available at
  \url{https://www.gurobi.com/wp-content/plugins/hd_documentations/documentation/9.0/refman.pdf}
  accessed on Jun. 14, 2021, 2020.

\bibitem{subramanyam2021joint}
Subramanyam, A., T.~Cokyasar, J.~Larson, and M.~Stinson.
\newblock Joint Routing of Conventional and Range-Extended Electric Vehicles in
  a Large Metropolitan Network, 2022.
\newblock Available at \url{https://arxiv.org/abs/2112.12769} accessed on Aug.
  8, 2022.

\end{thebibliography}

\end{document}